\documentclass[11pt]{article}

\usepackage{verbatim}
\usepackage{amsfonts,amsmath,amsthm}

\def\LaTeX{L\kern -.36em\raise .3ex\hbox{\sc a}\kern -.15em T\kern -.1667em%
\lower .7ex\hbox\mathbb{E}\kern -.125em X}

\usepackage{amssymb}
\usepackage{float}
\usepackage{epsfig}
\usepackage{amsmath}
\usepackage[english]{babel}
\usepackage{amstext}
\usepackage{mathrsfs}
\usepackage{amsthm}
\usepackage{fancyhdr}
\usepackage{latexsym}
\usepackage{amsmath}
\usepackage{amsfonts}
\usepackage{amssymb}
\usepackage{fancyhdr}

\usepackage[usenames]{color}
\definecolor{red}{rgb}{1.0,0.0,0.0}

\definecolor{blu}{rgb}{0.0,0.0,1.0}

\def\F{{\cal F}}

\def\ud{\mathrm{d}}

\def\to{\longrightarrow}

\def \0{{\textbf{0}}}
\def\norm{{\| \kern -.05em | }}

\newtheorem{Theorem}{Theorem}[part]
\newtheorem{Definition}{Definition}[part]
\newtheorem{Proposition}{Proposition}[part]
\newtheorem{Assumption}{Assumption}[part]
\newtheorem{Lemma}{Lemma}[part]
\newtheorem{Corollary}{Corollary}[part]
\newtheorem{Remark}{Remark}[part]

\def \N{\mathbb{N}}
\def \R{\mathbb{R}}

\def \E{\mathbb{E}}
\def \F{\mathbb{F}}

\def \P{\mathbb{P}}

\def \Ac{{\cal A}}

\def \Cc{{\cal C}}
\def \Dc{{\cal D}}

\def \Ic{{\cal I}}

\def \Lc{{\cal L}}

\def \Oc{{\cal O}}
 
\def \Sc{{\cal S}}
\def \Tc{{\cal T}}

\def \Yc{{\cal Y}}

\def \ep{\hbox{ }\hfill$\Box$}

\def\reff#1{{\rm(\ref{#1})}}

\def\beqs{\begin{eqnarray*}}
\def\enqs{\end{eqnarray*}}
\def\beq{\begin{eqnarray}}
\def\enq{\end{eqnarray}}

\begin{document}

\title{Characterization of the optimal boundaries in reversible investment problems}

\author{Salvatore FEDERICO\thanks{Part of this research was done when this author was Post-doc at  LPMA, Universit\'e Paris Diderot, supported by Alma Research.} 
\\\small  Dipartimento di Economia, 
\\\small Management  e Metodi Quantitativi
\\\small  Universit\`a di Milano 
\\\small  salvatore.federico at unimi.it
\and
Huy\^en PHAM
\\\small  Laboratoire de Probabilit\'es et
 \\\small  Mod\`eles Al\'eatoires, CNRS, UMR 7599
 \\\small  Universit\'e Paris Diderot
 \\\small  pham at math.univ-paris-diderot.fr
\\\small  and CREST-ENSAE
}

\date{This version:  July  2013}

\maketitle

\begin{abstract}
This paper studies a {\it reversible} investment problem where a social planner aims to control its capacity production in order to fit optimally the random demand of a good.  Our model allows for general diffusion dynamics on the demand as well as general cost functional. The resulting optimization problem leads  to a degenerate two-dimensional bounded variation  singular stochastic control problem, for which explicit solution is not available in general and the standard verification approach can not be applied a priori.  We use a  direct viscosity solutions approach for deriving some features of the optimal free boundary function, and for displaying  the structure of the solution. In the quadratic cost case, we are able to prove a smooth-fit $C^2$ property, which gives rise to a full characterization of the optimal boundaries and value function. 
\end{abstract}

\vspace{2cm}


\bigskip




\newpage


\section{Introduction}
We are concerned with a bounded variation singular  control problem motivated by a model of reversible investment.  More precisely, we imagine to deal with a social planner whose objective is to optimize some functional depending on the current demand of a good (energy, electricity, oil, corn, etc) and  its supply in terms of  production capacity that can be increased or decreased at any time and at  given  proportional costs.

Problems of  investment  under uncertainty have been introduced in the economic literature by \cite{M}  and then developed by several other authors (see \cite[Ch.\,11]{DP} for references on this subject). From a mathematical point of view, such  problems have been formulated as optimal stopping problems or, at a second stage of complexity,  as singular stochastic optimal control problems, and have given a considerable impulse to the development of the corresponding mathematical theory. As references for the theory of singular stochastic control in context different from  investment under uncertainty, we may mention the works  \cite{CMR, HS1, HS2, K2} and \cite[Ch.\,VIII]{FS}.  The mathematical literature of singular stochastic control applied to the subject of \emph{irreversible} investment under uncertainty (i.e. when the capacity can be only increased and the control is therefore monotone)  includes the works 
\cite{A1,  BK, B, CH1, CH2, DDSV, O, RS, W}. In particular \cite{B,RS} solve the problem by using a probabilistic representation result stated in \cite{BEL}, which seems very suitable to tackle this kind of problems, while \cite{W} uses a dynamic programming approach. The economic issue of  \emph{reversibility} (i.e. when the capacity can be also decreased and the control is a finite variation process) has then been introduced  and  studied, among others,  in \cite{AE, A2, GP, GT2,  LZ, MZ}. 
In the papers dealing with reversibility mentioned above, the ones (substantially) considering  two state variables (an uncontrolled one containing the noise, and a controlled one, representing the capacity) are \cite{A2, GT2, LZ, MZ}\footnote{We should mention also \cite{KW}, which just shows the connection between finite-variation singular control and Dynkin games. We shall indeed use this connection in Subsection \ref{sub:dynkin} to prove some results on the value function.}. 
 \cite{A2} derives optimality conditions based on economic considerations, while \cite{GT2} states and solves the problem with an interesting connection between finite-variation singular control problems and optimal switching problems. 
The papers dealing with  a dynamic programming approach  directly on the singular control problem  and  with the study of the associated Hamilton-Jacobi-Bellman equation (which in this case is a variational inequality) are \cite{LZ, MZ}. In particular, 
  \cite{MZ} considers an expected performance on infinite horizon with discounting over time, as in our case.  
However, the approach of \cite{MZ} is of \emph{verification} type. In a singular stochastic control framework, this  means that one has to guess some smooth fit properties of the value function at the optimal free boundary in order to look for a solution of the Hamilton-Jacobi-Bellman equation. Then one needs to prove \emph{a posteriori} that the solution found is indeed the value function and, as a byproduct, one gets also the optimal feedback control. When this approach is applicable, it turns out to be  very convenient, as it is theoretically fast (even if it may involve a nontrivial technical complexity) and allows a first understanding of the problem. Moreover, the presence of an explicit solution is an important tool to analyze the qualitative properties of  optimal control and trajectory. On the other hand, one has to recognize that it presents two drawbacks. First, it is based on a guess, and so it cannot bring to a deep understanding of the structural issues of the problem. Second,  it works only when explicit solutions are available, therefore it leaves the problem  completely unsolved most of the cases.

In the present paper, we perform a direct study of  the singular stochastic control  problem with bounded variation controls  (without passing through verification type arguments) by means of a viscosity approach to the Hamilton-Jacobi-Bellman (HJB) equation. To our knowledge, this is the first time that such an approach is used in the case of two state variables, in particular when the controlled state variable, here the reversible capacity process, has no diffusion term, and so is  degenerate\footnote{There are of course several papers (among them we may quote \cite{HS2}), which consider  singular stochastic control problems with multidimensional state variables, and characterize the value function in terms of viscosity solutions to the associated  HJB equations. However, rather few go beyond the viscosity characterization,  and investigate smooth-fit properties in order  to derive the structural form of the value function. 
In this spirit, we may mention the paper \cite{GP} in the case of just  one dimensional controlled variable.   See also \cite{GW} for impulse control of multi-dimensional diffusion processes with non degenerate diffusion  term. On the other hand,  we may quote the paper  \cite{SS}, which studies regularity  of a two-dimensional singular control problem with nondegenerate diffusion. Finally, we should mention the paper \cite{SS2}, dealing with a singular control problem with two state variables in a different context (consumption-investment under transaction costs). In this case the problem is approached by dynamic programming and by means of viscosity solutions to the associated Hamilton-Jacobi-Bellman equation. However, the regularity of the value function is proved by reducing the problem to dimension one, which is possible in that case due to the specific structure of the problem.}.  
Our approach allows us to keep much more generality with regard to the uncontrolled state variable (which is indeed a very general diffusion in the present paper, as in \cite{A2}) and to state the smooth-fit conditions of \cite{MZ} as \emph{necessary} conditions of optimality, i.e. prove that the value function \emph{must} satisfy these conditions\footnote{Another major advantage of such approach is that it allows generalizations. With this regard, we notice that here we minimize a cost functional. However, the arguments used here can be extended  to the case of profit/cost functional, as in \cite{MZ}.}. More precisely, we show that the value function is $C^1$ along the component of the controlled variable (Proposition \ref{prop:conv}; this easily follows from our assumptions  by convexity arguments, just working on the definition of value function). This allows to state the structure of the solution (Theorems \ref{prop:structure} and \ref{Verification}). Then, we prove that it has continuous mixed second derivative along the optimal boundary function (Proposition \ref{prop:sosf}; this is a deeper result, which invokes the viscosity property of the value function and requires the additional assumption \eqref{ggg} of quadratic cost in the capacity). The set of optimality conditions stated is then rewritten, following the arguments of \cite{A2}, in a more suitable way, which allows to determine the optimal boundaries, splitting them in three different regions and giving optimality conditions characterizing them in each of these regions (Theorem \ref{sol}). At the end, this machinery allows us to uniquely individuate the value function and solve the problem by Theorem \ref{Verification}. 
We mention that the approach developed in \cite{B} for singular control problem with monotone controls  is not valid anymore in the context of reversible investment.

The rest of the paper is organized as follows. In Section 2, we formulate  the two-dimensional bounded variation singular stochastic control problem and state the main assumptions. We study in Section 3 some first properties of the value function and of the optimal boundary, which is a function of the demand.  In Section 4, by relying on the viscosity property of the value function to its 
dynamic programming variational inequality, we  give a first main result providing the structure of the value function, and state a second main result  yielding the optimal control in terms of the optimal boundary. Section 5  focus on the case of quadratic  cost function, which allows us to prove a second order smooth fit principle. This  leads to the missing information to explicitly individuate the value function and the optimal boundary (the third main result), and makes the results of Section 4 applicable. Finally, we close the paper by explicit 
illustrations of the theory to the basic  example of  geometric Brownian motion for the uncontrolled demand diffusion in the case of irreversible investment. More examples and applications are developed, in the case of irreversible investment, in the companion paper \cite{AFPV}, where we also take into account delay in the expansion of the capacity production.

\section{The singular stochastic control problem}

\setcounter{equation}{0}
\setcounter{Assumption}{0} \setcounter{Theorem}{0}
\setcounter{Proposition}{0} \setcounter{Corollary}{0}
\setcounter{Lemma}{0} \setcounter{Definition}{0}
\setcounter{Remark}{0}

Let us fix  a probability space $(\Omega,\mathcal{F},\mathbb{P})$ equipped with a filtration $\F=(\mathcal{F}_t)_{t\geq 0}$ satisfying the usual conditions, and  supporting  
a standard one-dimensional Brownian motion   $(W_t)_{t\geq 0}$. 

On this space, we consider an uncontrolled state process  $D=(D_t)_{t\geq 0}$ (representing the demand of a good), governed by  a  diffusion dynamics:
\beq\label{eq:demand}
\ud  D_t &=& \mu (D_t) \ud  t + \sigma (D_t) \ud  W_t, \;\;\; D_0=d_0.
\enq
Let 
\beqs
\mathcal{O} &:=& (d_{\min},d_{\max}), \ \ \ -\infty\,\leq \, d_{\min}\,<\,d_{\max}\,\leq\,\infty.
\enqs
Throughout the paper we assume the following on the diffusion $D$.
\begin{Assumption}\label{ass:D1}
(i) The coefficients $\mu,\sigma:\Oc\rightarrow \R$ are  continuous and have at most linear growth.\\
(ii) For all $d_0\in \Oc$, there exists a unique non-exploding solution $D^{d_0}$ admitting a version with continuous path (and we shall always refer to such a version) to the SDE \eqref{eq:demand} in the space $(\Omega,\mathcal{F},\P)$ taking values into $\Oc$.\\
(iii) The unique solution $D$ continuously depends on the initial datum: if $d_n\stackrel{n\rightarrow \infty}{\to} d_0$, then $D^{d_n}\stackrel{n\rightarrow \infty}{\to} D^{d_0}$ almost surely.
\\
(iv) The SDE \eqref{eq:demand} satisfy a comparison criterion: if $d_0\leq d_0'$, then $D_t^{d_0}\leq D_t^{d_0'}$ $\P$-almost surely for every $t\geq 0$. \\
(v) The boundaries $d_{\min},d_{\max}$ are natural for the diffusion $D$ in the sense of Feller's classification and  the diffusion $D$ is regular.
\end{Assumption}
\begin{Remark}
{\rm Sufficient conditions for the assumptions above can be found in many classical references, such as, e.g.,  \cite[Ch.\,5]{KS}.}
{\rm We notice that some standard models of diffusion, such as arithmetic or geo\-metric Brownian motion, mean-reverting processes, or the CIR model (for suitable values of the parameters) satisfy Assumption \ref{ass:D1}.
}
\ep
\end{Remark}
Next,  we denote by $\mathcal{I}$ the class of  c\`adl\`ag bounded variation $\F$-adapted processes, setting $I_{0^-}$ $=$ $0$. 
Given $I\in\mathcal{I}$ we have the minimal decomposition $I=I^+-I^-$, where $I^+,I^-$ are the positive and the negative variation of $I$, respectively. It follows that the increments 
$$\ud I^+_t\ :=\ I^+_t-I^+_{t^-},  \ \ \ \ \ \ \ud I^-\ :=\ I^-_t-I_{t^-}$$ are supported on disjoint subsets of $[0,\infty)$.  We shall always refer to the latter minimal decomposition and, with a slight abuse of notation, we shall often denote $I=(I^+,I^-)$. 
The economic meaning of $I^+$ and $I^-$ is the following:
\begin{itemize}
\item[-]  $I^+_t$ is the cumulative investment done up to time $t$ to increase  the capacity;
\item[-] $I^-_t$ is the cumulative disinvestment done up to time $t$ to decrease the capacity.
\end{itemize}
Hence, the production capacity process  $(C_t)_{t\geq 0}$, controlled by $I$ $\in$ $\mathcal{I}$, is given by 
 \beq\label{GGG}
C_t &=& c_0 + {I_t^+}-{I_t^-}, \;\;\;\ \  c_0 \in \R.
\enq
The objective is to minimize over   $\mathcal{I}$
\beq \label{objective}
\mathbb{E}\Big[\int_0^{\infty} e^{-\rho t} \Big(g(C_t,D_t){{\ud}t}+{q_0^+}\ud {I_t^+}+{q_0^-}\ud {I_t^-}\Big)\Big],
\enq
 where $g$ $:$ $\mathbb{R}\times \mathcal{O}\rightarrow [0,\infty)$ is a cost function,  ${q_0^+}>0, \ q^-_0> 0$ are, respectively, the cost per unit of investment and the cost per unit of disinvestment, and $\rho$ is a positive discount factor.

\begin{Remark}\label{remrem}
{\rm
{\textbf{1.}} 
Among all the possible decompositions of a bounded variation process $I\in\mathcal{I}$, the minimal decomposition is the one providing the minimal value for the functional \eqref{objective}. 
{Indeed, denoting by $I^{m,+}-I^{m,-}$ the minimal decomposition of $I$, for all the other  decompositions $I=I^+-I^-$ the dynamics of the capacity $C$ is the same, while  $I^+\geq I^{m,+}, \ I^-\geq I^{m,-}$.} So 
\beqs
& &  \mathbb{E}\Big[\int_0^{\infty} e^{-\rho t} \Big(g(C_t,D_t){{\ud}t}+{q_0^+}\ud {I_t^{m,+}}+{q_0^-}\ud {I_t^{m,-}}\Big)\Big] \\
& \leq & \mathbb{E}\Big[\int_0^{\infty} e^{-\rho t} \Big(g(C_t,D_t){{\ud}t}+{q_0^+}\ud {I_t^+}+{q_0^-}\ud {I_t^-}\Big)\Big],
\enqs
\\ 
{\textbf{2.}} 
Even if we shall consider  $q_0^-$ as a finite number,  everything can be extended, giving a suitable sense,  to the case $q_0^-=\infty$. In this case the problem is  equivalent to require irreversibility for the investment (i.e. the case when $I^-$ is constrained to be $0$, as there is no convenience to disinvest, the cost being infinite). This  case is  treated in Subsection \ref{sec:irr}. 
\\ 
{\textbf{3.}}
For sake of simplicity, we do not impose the (economically meaningful: recall that $C$ should represent the capacity  production) state constraint $C_t\geq 0$. {We will comment in Remark \ref{Rem:positive} about the case that it may be verified a posteriori}. 
\\
\textbf{4.}
Note that, with respect to the usual  investment under uncertainty literature, which is mainly based on profit/cost  performance criterions, we focus here on the minimization of a  cost criterion in the spirit of a social planning problem, whose objective is to fit the capacity production to the demand at cheapest cost.
In particular the most significant case from the economic point of view is when {$g(c,d)=|c\,-\,d|^2$} (see also Remark \ref{rem:g}\,(2) below), as it represents a maximization of social surplus in the context of a linear inverse demand function (see \cite{AFPV} for a detailed description and explanation). We will give an explicit solution to the problem exactly in that case. \ep}
\end{Remark}

We shall make the following assumptions on the cost function $g$.

 \begin{Assumption}\label{ass:cost}
(i) {$g\in C^0(\R\times\Oc;\R_+)$, $g(\cdot,d)\in C^1(\mathbb{R};\R)$ for every $d\in\Oc$, and  $g_c\in C^0(\R\times \Oc; \R)$.} 

\noindent (ii)   $g(\cdot,d)$  is convex for all $d\in\mathcal{O}$ and $g_c(c,\cdot)$ is  {nonincreasing}  in $\mathcal{O}$ for every $c\in\mathbb{R}$. 

\noindent  (iii) $g$ and $g_c$ satisfy a polynomial growth condition w.r.t. $d$: there exist positive locally bounded functions $\gamma_0,\eta_0:\mathbb{R}\rightarrow \mathbb{R}$, and a constant $\nu\geq 0$ such that
\beq\label{growth:g}
|g(c,d)|+|g_c(c,d)| &\leq&  \gamma_0(c)+\eta_0(c)|d|^{\nu}, \ \ \ \forall \, c\in \mathbb{R}, \ \forall\,d\in {\mathcal{O}}.
\enq
\end{Assumption}

\begin{Remark}\label{rem:g}
{\rm 

\noindent {\bf 1.}  {We observe that the monotonicity assumption in Assumption \ref{ass:cost}-(ii) reflects an economic intuition.  
It means that the marginal cost with respect to capacity for a fixed level of capacity is nonincreasing in the demand: 
for a given level of capacity, the more is the demand, the more is convenient to invest; the less is the demand, the more is convenient to disinvest.}

\noindent {\bf 2.} Any function $g$ of the spread $|c-d|$ between capacity and demand,  in the form 
\beq\label{gspread}
g(c,d)&=& K_0 |c-d|^{\alpha}, \ \ \ K_0\geq 0, \ \alpha>1,
\enq
  satisfies Assumption \ref{ass:cost}.
\ep
}
\end{Remark}

\begin{Remark}\label{giorgio}
{\rm Following the idea of  \cite[Sec.\,6]{BK}, our model admits a suitable generalization to the case  of capacity dynamics in the form:
\beqs
\ud C_t &=& C_t( b \; \ud t + \gamma \; \ud W^0_t) + \ud I_{t}, \; \; C_{0^-} \; = \; c, 
\enqs
where $W^0$ is another Browinan motion independent of $W$. Indeed letting $C^0$ be the solution to  
\beqs
\ud C^0_t &=& C^0_t( b \; \ud t +\gamma \; \ud W^0_t), \;\; C^0_{0}\; = \; 1,
\enqs
the process $C$ can be rewritten as
\beqs
C_t \; = \; C^0_t \bar{C}_t,\ \ \ \ t\geq 0,
\enqs
where
\beqs
\bar{C}_t\ \ =\ \ c+\bar{I}^+_t-\bar{I}^-_t, & \mbox{ with } &  \bar{I}^+_t \; = \; \int_0^t\frac{1}{C^0_s}\,\ud I^+_s,\ \ \ \bar{I}^-_t \; = \; \int_0^t\frac{1}{C^0_s}\,\ud I^-_s.
\enqs
So, letting $\tilde{g}(\bar{c},c^0,d)= g(c^0\bar{c},d)$,  the problem becomes 
\beqs
\inf_{\bar{I}\in\mathcal{I}}\ \mathbb{E}\Big[\int_0^{\infty} e^{-\rho t} \Big(\tilde{g}(\bar{C}_t,C^0_t,D_t)\ud t+C_t^0({q_0^+}\ud \bar{I}^+_t+q_0^-\ud \bar{I}^-_t)\Big)\Big].
\enqs
This problem involves an additional uncontrolled state variable  (the variable $C^0$), but keeps the basic structures, so it seems approachable by the same techniques developed in the next sections. \ep
}
\end{Remark}

\section{Dynamic programming: preliminary results}

\setcounter{equation}{0}
\setcounter{Assumption}{0} \setcounter{Theorem}{0}
\setcounter{Proposition}{0} \setcounter{Corollary}{0}
\setcounter{Lemma}{0} \setcounter{Definition}{0}
\setcounter{Remark}{0}

We shall study the optimization problem by dynamic programming methods, and so we consider this singular stochastic control problem when varying initial data $(c_0,d_0)=(c,d)\in \R\times \mathcal{O}$. Therefore, from now on, we stress the dependence of  $C$ on $c,I$ and the dependence of $D$ on $d$ by denoting them respectively as  $C^{c,I}$, $D^d$. 
The state space is  then equal to
\beqs
\mathcal{S} &=& \mathbb{R}\times \mathcal{O}.
\enqs
Throughout the paper we indicate by $C^{h,k}(\mathcal{S};\R)$, $h, k\in\N$,  the class of functions which are  continuous, $h$-times differentiable with respect to the first variable, $k$-times differentiable with respect to the second variable, and having these derivatives  continuous in $\mathcal{S}$. 
\vspace{2mm}\\
Given $(c,d)\in\mathcal{S}$, the functional to be minimized over $I\in\mathcal{I}$ is
\beqs \label{functional}
G(c,d;I) &:=& \mathbb{E}\Big[\int_0^{\infty} e^{-\rho t} \Big(g(C^{c,I}_t,D^d_t){{\ud}t}+{q_0^+}\ud {I_t^+}+{q_0^-}\ud {I_t^-}\Big)\Big],
\enqs
and the associated value function  is
\beq\label{tildev}
{v}(c,d) &:=& \inf_{I\in\mathcal{I}}\ G(c,d;I), \;\;\; (c,d) \in \Sc. 
\enq
\subsection{First properties of the value function: finiteness and convexity}
Notice that  $v\geq 0$ as $g\geq 0$. We want to ensure also an upper bound for $v$.
Since $\mu,\sigma$ have at most linear growth, by  standard estimates we know (see, e.g.,  \cite[Ch.\,2.5, Cor.\,12]{K})  that there exist  constants {$K_{0}$ $=$ $K_{0,\mu,\sigma,\nu}\geq 0$ and  $K_1$ $=$ $K_{1,\mu,\sigma,\nu}$ $\in$ $\mathbb{R}$} such that
\beq\label{growth2}
\mathbb{E}\Big[ \big|D_t^d\big|^\nu \Big] &\leq& K_{0}(1+|d|^{\nu})e^{K_{1} t},  \;\;\; \forall  t\geq 0.
\enq
In the sequel, we make the standing assumption that the discount factor $\rho$ satisfies
\beq \label{rhoK1}
\rho & > & {K_1^+},
\enq
where  $K_{1}$ is the constant appearing in \eqref{growth2}. 
Using Assumption \ref{ass:cost}\,(iii) and \eqref{growth2}-\eqref{rhoK1}, we get 
 \beq \label{mingamma}
   \hat{V}(c,d) \ := \ \mathbb{E}\Big[ \int_0^{\infty} e^{-\rho t} g(c,D_t^d){\ud}t \Big]\ \leq \ \gamma_1(c)+\eta_1(c)|d|^\nu, \ \ \ \forall (c,d)\in\Sc,
\enq
for some nonnegative locally bounded real functions $\gamma_1,\eta_1$.
Moreover, due  Assumption \ref{ass:cost}, the function $\hat{V}$ is continuous in $\mathcal{S}$ and differentiable with respect to $c$ for all $d\in\Oc$, with
\beq\label{esthatV}
  \hat{V}_c(c,d) & =  & \mathbb{E}\Big[ \int_0^{\infty} e^{-\rho t} g_c(c,D_t^d){\ud}t \Big],  \ \ \ \ (c,d)\in\Sc,
\enq
and for the same reason as before
\beq\label{esthatV'}
\hat{V}_c(c,d) &\leq& \gamma_1(c)+\eta_1(c)|d|^\nu, \ \ \ \forall (c,d)\in\Sc.
\enq

Now,
let $d_0\in\Oc$ be a reference point and let us introduce the functions
$$
S'(d)\ \ :=\ \ \mbox{exp}\left(-\int_{d_0}^d \frac{2\mu (\xi)\ud\xi}{\sigma^2(\xi)}\right),\ \ \ \  d\in\Oc,
$$
and 
$$
m'(d)\ \ :=\ \  \frac{2}{\sigma^2(d)S'(d)}, \ \ \ \ d\in \Oc.
$$
$S'$ is the the density of the so called scale function of the diffusion $D$, and $m'$  is the density of the so called speed measure of the diffusion $D$.
Let us denote respectively by 
$\psi$ and $\varphi$ the increasing and   decreasing fundamental solutions, individuated up to a multiplicative constant, to the linear ordinary differential equation
\beq\label{ODE}
\Lc \phi(d) \  := \   \rho \phi(d)-\mu(d)\phi'(d)-\frac{1}{2}\sigma^2(d)\phi''(d) &=& 0.
\enq
The existence and properties of such 
functions, as well as their relationship with the functions $S, m$ defined above, can be found in several references including in  \cite[Ch.\,II]{BS}, \cite[Ch.\,15]{KT}, \cite[Ch.\,V]{Rog}, and \cite[Ch.\,2]{Mandl}. 
In particular we know that $\psi, \varphi$  are strictly positive, convex, and,  since $d_{\min}, d_{\max}$ are natural boundaries, they satisfy (see, e.g., \cite[Ch.\,2]{BS})
\beq \label{psiphi}
\lim_{d\downarrow d_{\min}}\psi(d) \; = \; 0, \ \ \  \lim_{d\downarrow d_{\min}}\varphi(d) \; = \; \infty, \ \ \  \lim_{d\uparrow d_{\max}}\psi(d) \; = \; \infty, \ \ \   \lim_{d\uparrow d_{\max}}\varphi(d) \; = \; 0,
\enq
\beq \label{psiphi'}
\lim_{d\downarrow d_{\min}}\frac{\psi'(d)}{S'(d)} \; = \; 0, \ \ \  \lim_{d\downarrow d_{\min}}\frac{\varphi'(d)}{S'(d)} \; = \; -\infty, \ \ \  \lim_{d\uparrow d_{\max}}\frac{\psi'(d)}{S'(d)} \; = \; \infty, \ \ \   \lim_{d\uparrow d_{\max}}\frac{\varphi'(d)}{S'(d)} \; = \; 0.
\enq
Let $w$ be the constant positive Wronskian of the fundamental solutions $\psi,\varphi$, i.e. 
$$0\ \ <\ \ w\ \ \equiv \ \ \frac{\psi'(d)\varphi(d)-\psi(d)\varphi'(d)}{S'(d)},\ \ \ \ d\in \Oc.$$
and let $p(t,d,\cdot)$ be the density of the transition probability $P(t,d,\cdot)$ of the diffusion $D$.
Using the characterization of the Green's function 
$$G(d,h) \ \ :=\ \ \int_{0}^\infty e^{-\rho t} p(t,d,h)\ud t$$
associated to $D$ as
$$
G(d,h)\ \ =\ \ \begin{cases}
w^{-1} \psi(d)\varphi(h), \ \ \ \ \mbox{if} \ d\leq h,\\
w^{-1} \psi(h)\varphi(d), \ \ \ \ \mbox{if} \ d\geq h,
\end{cases}
$$
 and the fact that it is the kernel of the resolvent operator (see, e.g., \cite[Ch.\,V]{Rog} or \cite[Ch.\,15]{KT}) with respect to $m$, i.e. 
$$\mathbb{E}\left[\int_0^\infty e^{-\rho t} f(D_t^d)\ud t\right] \ \ =\ \ \int_\mathcal{O} f(h)G(d,h)m'(h)\ud h, \ \ \ \ \forall f\in \mathcal{B}(\mathcal{O};\R),$$
 we see (approximating $g,g_c$  by bounded functions and using the monotone convergence theorem) that the functions $\hat{V}$ and $\hat{V}_c$ can be represented in terms of $\psi,\varphi$ as
\beq\label{Vc}
\hat{V}(c,d) & =& w^{-1}\Big[\varphi (d)\int_{d_{\min}}^{d} \psi(\xi)g(c,\xi)m'(\xi)\ud \xi\ +\  \psi (d)\int_{d}^{d_{\max}} \varphi(\xi)g(c,\xi)m'(\xi)\ud \xi\Big],\\
\hat{V}_c(c,d) & =& w^{-1}\Big[\varphi (d)\int_{d_{\min}}^{d} \psi(\xi)g_c(c,\xi)m'(\xi)\ud \xi\ +\  \psi (d)\int_{d}^{d_{\max}} \varphi(\xi)g_c(c,\xi)m'(\xi)\ud \xi\Big],\label{Vcd}
\enq

\begin{Proposition}\label{prop:conv}
The value function $v$ is convex with respect to $c$ and satisfies the growth condition, for some locally bounded functions $\gamma_1,\eta_1:\R\to\R$,
\beq \label{growthv}
0  \ \, \leq\ \, v(c,d) & \leq & \hat{V}(c,d) \ \ \leq\ \ \gamma_1(c)+\eta_1(c) |d|^{\nu}, \ \ \ \forall (c,d)\in \mathcal{S},
\enq
\end{Proposition}
\noindent \textbf{Proof.} 
\eqref{growthv} comes from \eqref{esthatV} and from the inequality $v(c,d)\leq G(c,d;0)=\hat{V}(c,d)$.

Convexity of $v$  follows in a standard way from the convexity of $g$ with respect to $c$ and linearity of the state equation for $C^{c,I}$. \ep

%
%
\subsection{Existence of optimal controls and the associated Dynkin game}\label{sub:dynkin}

In this subsection we show that the singular stochastic control problem admits optimal controls and that it is related to a suitable associated Dynkin game. 
We estabilish this connection mainly to  inherit from the monotonicity of $g_c(c,\cdot)$ the monotonicity of  $v_c(c,\cdot)$, whose direct proof seems not attainable. 
The proofs of Propositions \ref{exist}, \ref{Prop:dynkin} closely follow the arguments of \cite{KW}, and are reported in Appendix. 

\begin{Definition}
Given $(c,d)\in\mathcal{S}$ we say that a control $I^*\in\mathcal{I}$ is optimal starting from $(c,d)$ if $G(c,d;I^*)=v(c,d)$.
\end{Definition}
\begin{Proposition}\label{exist}
For all $(c,d)\in\mathcal{S}$ there exists an optimal control $I^*$ starting from $(c,d)$. Moreover, if $g(\cdot,d)$ is strictly convex on $\mathbb{R}$ for every $d\in\Oc$, then $I^*$ is the unique (up to undistinguishability) optimal control starting from $(c,d)$. 
\end{Proposition}
\
Let $\mathcal{T}$ denote the set of all $\mathbb{F}$-stopping times.
For fixed  $(c,d)\in\mathcal{S}$, we  may consider the functional, controlled by   $\sigma\in\mathcal{T},\ \tau\in\mathcal{T}$,
\beq\label{fun:game}
J(c,d;\sigma,\tau)&=&\E\Big[\int_0^{\sigma\wedge\tau} e^{-\rho t} g_c(c,D_t^d)\ud t+q_0^-e^{-\rho \sigma}\mathbf{1}_{\{\sigma<\tau\}}-q_0^+e^{-\rho \tau}\mathbf{1}_{\{\tau<\sigma\}}\Big].
\enq
We can imagine that $J(c,d;\tau,\sigma)$ is the payoff associated to a two-players stochastic game. The two players, P1 and P2, have the possibility to stop the game at times $\sigma$ and $\tau$, respectively (i.e. P1 controls the game through $\sigma$ and P2 controls the game through $\tau$). If P1 stops first ($\sigma<\tau$), he pays to P2 the amount $q_0^-e^{-\rho \sigma}$; if P2 stops first ($\tau<\sigma$), he pays to P1 the amount $q_0^+e^{-\rho \tau}$; if they decide to stop at the same time, i.e. if $\tau=\sigma$, then no cashflow occurs; finally, as long as the game is running, i.e. up to time $\sigma\wedge\tau$, P1 pays P2 at the rate $e^{-\rho t}g_c(c,D_t^d)$ per unit of time. The goal of P1 is to minimize \eqref{fun:game}, while the goal of P2 is to maximize \eqref{fun:game}. The functions
\beqs
\underline{w}(c,d)\ \ :=\ \ \sup_{\tau\in\mathcal{T}}\,\inf_{\sigma\in\mathcal{T}} \ J(c,d;\sigma,\tau), \ \ \ \ \ \ \ \ \overline{w}(c,d)\ \ :=\ \ \inf_{\sigma\in\mathcal{T}}\,\sup_{\tau\in\mathcal{T}}\ J(c,d;\sigma,\tau),
\enqs
are called \emph{lower- and upper-values of the game}.  Clearly $\underline{w}(c,d)\leq \overline{w}(c,d)$. If $\underline{w}(c,d)=\overline{w}(c,d)$, the game is said to have a \emph{value} denoted by $w(c,d):=\underline{w}(c,d)=\overline{w}(c,d)$. A pair $(\sigma^*,\tau^*)\in\mathcal{T}\times\mathcal{T}$ is called a \emph{saddle-point of the game} if
\beq\label{saddle}
J(c,d;\sigma^*,\tau)\ \ \leq \ \ J(c,d;\sigma^*,\tau^*)\ \ \leq \ \  J(c,d;\sigma,\tau^*), \ \ \ \ \  \ \forall \sigma\in\mathcal{T}, \ \ \forall \tau\in\mathcal{T}.
\enq
One easily sees that the existence of a saddle point implies that the game has a value and
\beq
w(c,d)&=& J(c,d;\sigma^*,\tau^*).
\enq
\begin{Proposition}\label{Prop:dynkin}
\begin{enumerate}\item
Let $(c,d)\in\mathcal{S}$ and let $I^*=(I^{*,+}, I^{*,-})\in\mathcal{I}$ be an optimal control for the singular stochastic control problem, i.e. such that $v(c,d)=G(c,d;I^*)$. Define the stopping times 
\beqs
\sigma^*\ \ :=\ \ \inf\,\{\,t\geq 0 \ \ | \ \ I^{*,-}_t>0\,\}, \ \ \ \  \ \ \tau^*\ \ :=\ \ \inf\, \{\,t\geq 0 \ \ |\  \ I^{*,+}_t>0\,\}.
\enqs
Then $(\sigma^*,\tau^*)\in\mathcal{T}\times\mathcal{T}$ is a saddle point for the associated Dynkin game.
\item $v$ is differentiable with respect to $c$ in $\mathcal{S}$ and it holds the equality $v_c=w$, where $w$ is the (well-defined) value of the associated Dynkin game.
\end{enumerate}
\end{Proposition}

\vspace{2mm}

By relying on  this connection between singular control and Dynkin game, we prove now some pro\-perties on the derivative of the value function $v_c$, to be used in the next Section.

\begin{Proposition}\label{prop:vc}
The function $v_c$ has the  following properties:
\begin{enumerate}
\item
$v_c$ is continuous in $\mathcal{S}$.
\item
$v_c(c,\cdot)$ is nonincreasing in $\mathcal{O}$ for all $c\in\R$.
\item
$-q_0^+\ \leq \ v_c\ \leq \ {q_0^-}$ in $\mathcal{S}$.
\end{enumerate}
\end{Proposition}
\noindent \textbf{Proof.}   
 1. Let $(c,d)\in\mathcal{S}$ and take a sequence $(c_n,d_n)\rightarrow (c,d)$. For each $n\in\N$, let $(\sigma^{*}_n,\tau^*_n)$ be a saddle-point for the Dynkin game starting at $(c_n,d_n)$, and let $(\sigma^*,\tau^*)$ be a saddle point for the Dynkin game starting at $(c,d)$. Using \eqref{saddle}, we then have
 \beq\label{pasd}
 w(c,d)-w(c_n,d_n)&=& J(c,d;\sigma^*,\tau^*)-J(c_n,d_n;\sigma^*_n,\tau^*_n)\nonumber\\
  &\leq& J(c,d;\sigma^*_n,\tau^*)-J(c_n,d_n;\sigma^*_n,\tau^*)\nonumber \\
  &=&\E\left[\int_0^{\tau^*\wedge\sigma^*_n} e^{-\rho t} \big(g_c(c,D_t^d)-g_c(c_n,D_t^{d_n})\big)\ud t\right]\\
  &=&\E\left[\int_0^{\infty} e^{-\rho t} \big(g_c(c,D_t^d)-g_c(c_n,D_t^{d_n})\big)\mathbf{1}_{\{t\leq \tau^*\wedge\sigma^*_n\}}\ud t\nonumber
  \right].
 \enq
Note that, assuming without loss of generality that $(d_n)_{n\in\N}\subset (d-\varepsilon,d-\varepsilon)\subset\Oc$ for suitable $\varepsilon>0$, we have by Assumption \ref{ass:D1}\,(iv) 
\beq   \label{qqa}
|D_t^{d_n}|&\leq& |D_t^{d-\varepsilon}|+|D_t^{d+\varepsilon}|, \ \ \ \forall t\geq 0, \ \forall n\in\N.
\enq 
On the other hand Assumption \ref{ass:D1}\,(iii)  ensures the convergence
\beq\label{qqw}
D_t^{d_n} \ \stackrel{n\rightarrow\infty}{\to} \ D_t^d, \ \ \ \ \mbox{a.s.},\  \ \forall t\geq 0.
\enq
Hence, using Assumption \ref{ass:cost}, \eqref{rhoK1}, and \eqref{qqa}-\eqref{qqw},  we can  apply dominated convergence to \eqref{pasd} for $n\rightarrow\infty$ and conclude that $\liminf_{n\rightarrow\infty} w(c_n,d_n)\geq w(c,d)$. 

Arguing in a similar way, but considering the couple $(\sigma^*,\tau_n^*)$ in place of the couple $(\sigma_n^*,\tau^*)$, one also gets the  inequality $\limsup_{n\rightarrow\infty} w(c_n,d_n)\leq w(c,d)$, so $w$ is continuous at $(c,d)$. 

Then the claim follows by Proposition \ref{Prop:dynkin}\,(2). 
 \smallskip
 
 2. 
By the assumption that $g_c(c,\cdot)$ is  {non increasing} (Assumption \ref{ass:cost}(ii)), 
and from the same comparison result cited above, we have, for every $d,d'\in\mathcal{O}$ such that $d\leq d'$, 
$$
J(c,d;\sigma,\tau)\ \ \geq \ \ J(c,d';\sigma,\tau), \ \ \ \ \ \ \forall \sigma\in\mathcal{T}, \ \forall \tau\in\mathcal{T}.
$$
Passing to the infimum over $\sigma\in \mathcal{T}$ and then to the supremum over $\tau\in \mathcal{T}$ the inequality above we get, for every $d,d'\in\mathcal{O}$ such that $d\leq d'$,
$$
\underline{w}(c,d)\ \ \geq \ \ \underline{w}(c,d').
$$
Proposition \ref{Prop:dynkin} states that the game has a value, so from the inequality above we get,
for every $d,d'\in\mathcal{O}$ such that $d\leq d'$,
$$
{w}(c,d)\ \ \geq \ \ {w}(c,d').
$$
Hence, the claim follows from Proposition \ref{Prop:dynkin},(2). 
 \smallskip
 
 3. We have  $J(c,d;\sigma,0)=-q_0^+$  for every $\sigma\in\mathcal{T}$, and $J(c,d;0,\tau)=q_0^-$ for every $\tau\in\mathcal{T}$. It follows that $-q_0^+\leq w(c,d)\leq q_0^-$ and the claim follows from Proposition \ref{Prop:dynkin}\,(2).
%
%
\ep

\section{The dynamic programming equation and the structure of the solution}

\setcounter{equation}{0}
\setcounter{Assumption}{0} \setcounter{Theorem}{0}
\setcounter{Proposition}{0} \setcounter{Corollary}{0}
\setcounter{Lemma}{0} \setcounter{Definition}{0}
\setcounter{Remark}{0}

In view of Proposition \ref{prop:vc}, we introduce  the so-called \emph{continuation region}
\beqs
\mathcal{C} &:=& \{\, (c,d)\in {\mathcal{S}} \ \,  | \ -q_0^+<\,v_c(c,d)\, <\, q_0^-\,\},
\enqs
and its complement set,  the \emph{action region}
\beq\label{mathD}
\mathcal{A} &:=& \mathcal{A}^+ \ \cup \ \mathcal{A}^-,
\enq
where $\mathcal{A}^+$ and $\mathcal{A}^-$ are respectively the \emph{investment} and the \emph{disinvestment region} defined by
\beq\label{A+-}
\mathcal{A}^+ \  := \   \{(c,d)\in{\mathcal{S}} \ | \ v_c(c,d)=-{q_0^+}\}\,, \ \ \ \mathcal{A}^- \ := \   \{(c,d)\in{\mathcal{S}} \ | \ v_c(c,d)={q_0^-}\}.
\enq
We also set
$$
\partial^+\Cc \ \:=\ \ \bar{\Cc}\cap \mathcal{A}^+, \ \ \ \ \partial^-\Cc \ \:=\ \ \bar{\Cc}\cap \mathcal{A}^-.
$$ 
The boundaries  $\partial^{\pm}\Cc$  are associated with a free boundary differential problem (which we are going to define in the next subsection) and are the objects to individuate to solve the optimal stochastic control problem.\\\\
Let us then consider  the functions $\hat c_+, \hat c_-$ $:$ $\Oc$ $\rightarrow$ $\bar\R$ defined with the conventions $\inf\emptyset=\infty$, $\inf \mathbb{R}=-\infty$, $\sup \R =\infty$, $\sup\emptyset=-\infty$ (the equalities below are consequence of convexity of $v$ with respect to $c$):
\beq\label{c+}
{\hat{c}_+}(d) & := &   \inf\ \{ c\in \mathbb{R} \  | \ v_c(c,d)>-{q_0^+}\}\ \ =\ \  \sup\ \{ c\in \mathbb{R} \  | \ v_c(c,d)=-{q_0^+}\},\\ 
 {\hat{c}_-}(d) &  :=&   \sup\ \{ c\in \mathbb{R} \  | \ v_c(c,d)< {q_0^-}\}\ \ =\ \  \inf\ \{ c\in \mathbb{R} \  | \ v_c(c,d)={q_0^-}\}.\label{c-}
\enq

%
%
\begin{Proposition}\label{frontiera}
\begin{enumerate}
\item $\hat{c}_+: \Oc\rightarrow \R\cup\{-\infty\}, \ \hat{c}_-: \Oc\rightarrow \R\cup\{\infty\}$, they are both nondecreasing and
\beq\label{c+-}
\hat{c}_+(d) \ \ <\ \ \hat c_-(d), \ \  \ \ \ \ \forall d\in\mathcal{O}.
\enq
\item  $\hat{c}_+$ is right-continuous and $\hat{c}_-$ is left-continuous. 
\item
The action and continuation regions are expressed in terms 
of the functions $\hat c_{\pm}$ as: 
\beqs\label{mathC1}
\mathcal{C} & = & \{(c,d)\in {\mathcal{S}} \ \ | \ \ \hat{c}_+(d) \ <\  c\ <\ {\hat{c}_-}(d)\},
\enqs
\vspace{-1.4cm}\\
\beqs 
\mathcal{A}^+ &  = &  \{(c,d)\in{\mathcal{S}} \ \ |\  \ c \ \leq \ {\hat{c}_+}(d)\}, \ \ \ \ \ 
\mathcal{A}^- \ \   = \ \  \{(c,d)\in{\mathcal{S}} \ \ |\ \ c\ \geq \ {\hat{c}_-}(d)\}.
\enqs
\item $\Cc$ is open and connected, and $\mathcal{A}^{\pm}$ are closed and connected.
\end{enumerate}
 \end{Proposition}
\noindent \textbf{Proof.} 
1. The fact that $\hat{c}_+$ takes values in $\R\cup\{-\infty\}$ and $\hat{c}_-$ takes values in $\R\cup\{\infty\}$ is consequence of the nonnegativity of $v$, combined with the convexity of $v(\cdot,d)$ and with \eqref{c+}-\eqref{c-}. 
Monotonicity follows from Proposition \ref{prop:vc}\,(2) and  \eqref{c+}-\eqref{c-}.  Finally, \eqref{c+-} is due to the convexity of $v$ with respect to $c$ and to the fact that $v(\cdot,d) \in C^1(\mathbb{R};\R)$ for every $d\in\Oc$. 

2. It follows from Proposition \ref{prop:vc}\,(1) and from the convexity of $v$ w.r.t. $c$. 

3-4. They follow from the previous items also considering  \eqref{c+}-\eqref{c-}. \ep
\vspace{2mm}

\noindent Below it is represented a possible shape of the regions $\Cc,\Ac^{\pm}$ and of the  functions $\hat{c}^\pm$ (here $d_{\max}=\infty$).

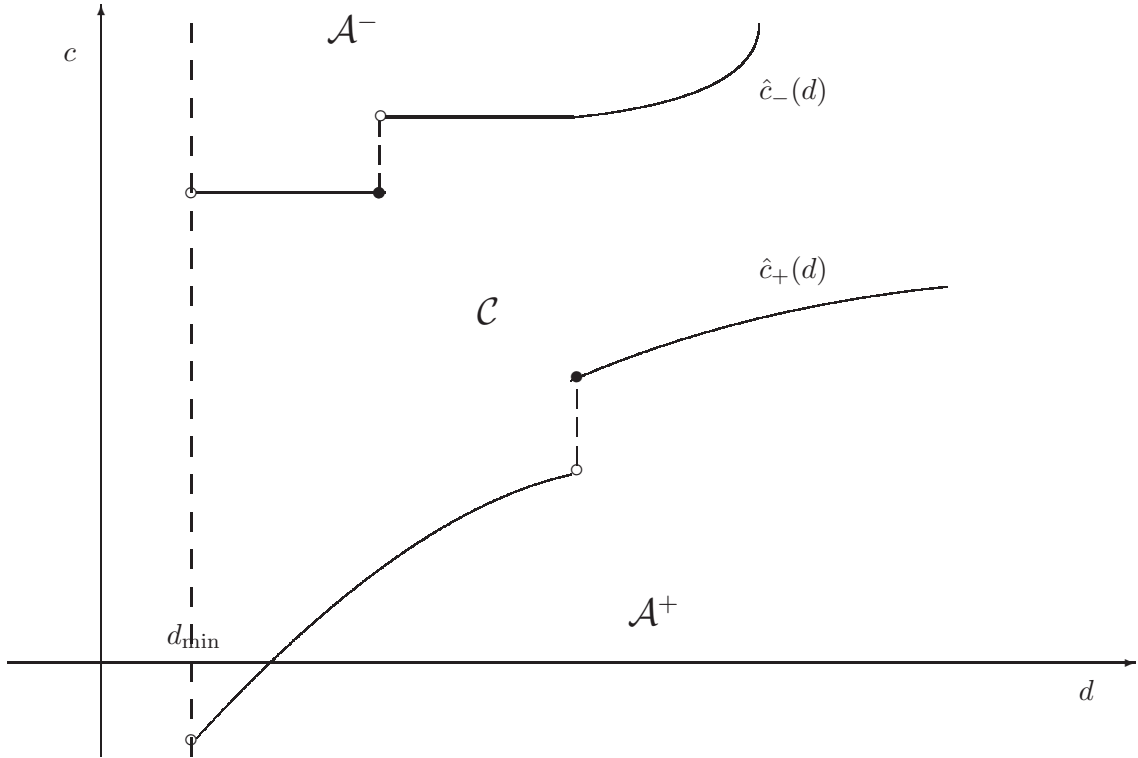
 \setlength{\unitlength}{2.5cm}
\begin{picture}(6,5)(-3,-2)
\put(-3,-1){\vector(1,0){6}}
\put(-1.99,1.5){\line(1,0){1}}
\put(-1.99,1.501){\line(1,0){1}}
\put(-1.99,1.502){\line(1,0){1}}
\put(2.7,-1.2){$d$}
\put(-2.06,1.46){$\circ$}
\put(-1.06,1.46){$\bullet$}
\put(-.99,1.9){\line(1,0){1}}
\put(-.99,1.901){\line(1,0){1}}
\put(-.99,1.902){\line(1,0){1}}
\put(-.99,1.906){\line(1,0){1}}
\put(-1.05,1.87){$\circ$}
\put(-1.06,1.46){$\bullet$}
\put(-2.15,-0.9){$d_{\min}$}
\put(-.5,.8){$\mbox{\Large{$\mathcal{C}$}}$}
\put(.3,-.8){$\mbox{\Large{$\mathcal{A}^+$}}$}
\put(-1.3,2.3){$\mbox{\Large{$\mathcal{A}^-$}}$}
\put(-.01,-.012){$\circ$}
\put(-2.06,-1.45){$\circ$}
\put(-.01,.48){$\bullet$}
\put(-2.5,-1.5){\vector(0,1){4}}
\put(1,1.05){${\hat{c}_+}(d)$}
\put(1,2){${\hat{c}_-}(d)$}
\put(0.03,0.05){\line(0,1){.1}}
\put(0.03,0.20){\line(0,1){.1}}
\put(0.03,0.35){\line(0,1){.1}}
\put(-1.02,1.50){\line(0,1){.1}}
\put(-1.02,1.65){\line(0,1){.1}}
\put(-1.02,1.80){\line(0,1){.08}}
\put(-2.7,2.2)
{$c$}
\multiput(-2.02,-1.5)(0,.2){20}{\line(0,1){.1}}
\qbezier(0,.5)(0.8853,0.8853)
(2,1)
\qbezier(0,0)(-0.8853,-.1953)
(-1.99,-1.4)
\qbezier(0,.5)(0.8853,0.8853)
(2,1)
\qbezier(0,0)(-0.8853,-.1953)
(-1.99,-1.4)
\qbezier(0,.5)(0.8853,0.8853)
(2,1)
\qbezier(0,0)(-0.8853,-.1953)
(-1.99,-1.4)
\qbezier(0.0,1.9)(1,2)
(1,2.4)
\qbezier(0.0,1.9)(1,2)
(1,2.4)
\qbezier(0.0,1.9)(1,2)
(1,2.4)\qbezier(0.0,1.9)(1,2)
(1,2.4)
\qbezier(0.0,1.9)(1,2)
(1,2.4)
\qbezier(0.0,1.9)(1,2)
(1,2.4)
\end{picture}

\vspace{2mm}


\noindent
Let us define
\beqs
\underline{c}_+ \; := \; \inf_{d \in \mathcal{O}} {\hat{c}_+}(d), \ \ \ \
\bar{c}_+ \; := \;  \sup_{d\in \mathcal{O}} {\hat{c}_+}(d), \ \ \ \underline{c}_- \; := \; \inf_{d \in \mathcal{O}} {\hat{c}_-}(d), \ \ \ 
\bar{c}_- \; := \;  \sup_{d\in \mathcal{O}} {\hat{c}_-}(d),
\enqs
and the pseudo-inverse  of  ${\hat{c}_{\pm}}$, i.e. the functions  $\hat{d}_{\pm}:\mathbb{R}\to\bar{\Oc}$ defined by 
\beq\label{pseudo}
\hat{d}_+(c) \ :=\ \inf\ \{d\in\Oc \ | \ \hat{c}_+(d)\geq  c\}, \ \ \ \  \hat{d}_-(c) \ :=\ \sup\ \{d\in\Oc \ | \ \hat{c}_-(d)\leq c\},
\enq
with the convention $\inf\emptyset=d_{\max}$ and $\sup\emptyset=d_{\min}$.  
\begin{Proposition}\label{frontierad}
\begin{enumerate}
\item We have the equalities
\beq\label{chard}
\hat{d}_+(c) \ =\  \sup\ \{d\in\mathcal{O} \ | \ v_c(c,d)>-{q_0^+}\},  \ \ \ \ \hat{d}_-(c)  \ = \ \inf\ \{d\in\mathcal{O} \ | \ v_c(c,d)<{q_0^-}\}.
\enq
\item 
The functions $\hat{d}_{\pm}$ are nondecreasing and $\hat{d}_+\geq \hat{d}_-$.
\item
If  $\bar{c}_-<\infty$, then $\hat{d}_{-}=d_{\max}$ on $[\bar{c}_-,\infty)$; if $\underline{c}_+>-\infty$, then $\hat{d}_{+}=d_{\min}$ on $(-\infty,\underline{c}_+]$.
\item
$\hat d_- (c)<  \hat d_+(c)$ if and only if $c \in (\underline{c}_+, \bar{c}_-)$.
\end{enumerate}
\end{Proposition}
\noindent
\textbf{Proof.}
1. It directly follows from 
 the definition of $\hat{c}_\pm$, $\hat{d}_{\pm}$.
 
 2. Monotonicity of  $\hat{d}_\pm$ and the inequality $\hat{d}_+\geq \hat{d}_-$ follow  from Proposition \ref{frontiera}\,(1).  
 
 3. By monotonicity of $\hat{d}_-$,   $\lim_{c\rightarrow\infty+}\hat{d}_-(c)$ exists. 
Suppose by contradiction $\lim_{c\rightarrow\infty}\hat{d}_-(c)=\bar{d}<d_{\max}$. This would imply $\hat{c}_- =\infty$ over $(\bar{d},d_{\max})$, which contradicts $\bar{c}_-<\infty$. A similar argument works for the other claim.

4. It follows from 
  \eqref{c+-}.
\ep
 
\vspace{2mm}

We also introduce  the $c$-section sets of the continuation region
\beq\label{secc}
S_c &:= &\{c\}\times (\hat{d}_-(c),\hat{d}_+(c)), \ \ \ c \in\R.
\enq
Due to Proposition \ref{frontierad}, we have
\beq\label{snemp}
 c \in (\underline{c}_+, \bar{c}_-) \ \Longleftrightarrow  \ \hat{d}_-(c)<\hat{d}_+(c) \   \Longleftrightarrow \ S_c\neq \emptyset.
\enq
We have the following  result on the form of the continuation region.
\begin{Proposition}\label{sezioni}
 We have the representation of the continuation region
\beq\label{cupcirc}
{\mathcal{C}} &=& \bigcup_{c\in (\underline{c}_+,\bar{c}_-)}S_c.
\enq
\end{Proposition}
\noindent {\bf Proof.} 
If $(c,d)\in\ {\mathcal{C}}$, then $-q_0^+<v_c(c,d)<q_0^-$, so, by continuity of $v_c$ (Proposition \ref{prop:vc}\,(1)), it is  $-q_0^+< \hat{v}_c< {q_0^-}$ in some suitable neighborhood of $(c,d)$. Then $\hat{d}_-(c)<\hat{d}_+(c)$,  therefore, by \eqref{snemp},  $c\in  (\underline{c}_+,\bar{c}_-)$ and $(c,d)\in S_c\neq\emptyset$. Hence we have proved the inclusion ${\mathcal{C}}\ \subset\ \bigcup_{c\in (\underline{c}_+,\bar{c}_-)}S_c$.

Conversely, let $c\in (\underline{c}_+,\bar{c}_-)$ and let $d\in\Oc$ be such that $(c,d)\in S_c(\neq \emptyset)$. 
{By \eqref{chard} and \eqref{snemp}, we have $-q_0^+< v_c(c,\cdot)<{q_0^-}$ in some neighborhood of $d$. 
The  continuity of $v_c$ with respect to $c$ (Proposition \ref{prop:vc}\,(1)) implies   $-q_0^+< v_c<{q_0^-}$ in some neighborhood of $(c,d)$. 
Therefore  $(c,d)\in\Cc$. Hence we have proved the inclusion ${\mathcal{C}}\ \supset\ \bigcup_{c\in (\underline{c}_+,\bar{c}_-)}S_c$. }
\ep
\\\\
We also introduce the functions $\hat{c}_{\pm,g}$ from $\Oc$ into $\overline\R$ defined, with the usual convention $\sup \emptyset =-\infty, \inf\emptyset =\infty$, by:
\beqs
\hat{c}_{+,g}(d) \; = \; \inf\{c\in\mathbb{R} \ | \ g_c(c,d)>-\rho q_0^+\}, & &  
\hat{c}_{-,g}(d) \; = \; \sup\{c\in\mathbb{R} \ | \ {g}_c(c,d)<\rho q_0^-\}.
\enqs
One easily checks that, by Assumption \ref{ass:cost}, they are nondecreasing  and, respectively  right- and left-continuous. Moreover, we clearly have, by convexity of $g(\cdot,d)$ and continuity of $g_c$, the inequality $\hat{c}_{+,g}<\hat{c}_{-,g}$. We have the following estimates of $\hat{c}_{\pm}$ in terms of $\hat{c}_{\pm,g}$.
\begin{Proposition}\label{prop:stime}
$\hat{c}_+\leq \hat{c}_{+,g}$ and $\hat{c}_-\geq \hat{c}_{-,g}$.
\end{Proposition}
\noindent \textbf{Proof.}
Let us show the first inequality, the second one can be proved symmetrically.
Let
 $d\in{\mathcal{O}}$ and take $c> \hat{c}_{+,g}(d)$, so that $g_c(c,d)+\rho q_0^+>0$. Let $\varepsilon\in\left(0,\frac{g_c(c,d)+\rho q_0^+}{\rho}\right)$, and consider the stopping time
 $$\tau_\varepsilon \ \ := \ \ \inf\,\{t\geq 0\ | \ g_c(c,D_t^d)+\rho q_0^+\leq \rho\varepsilon\,\}.$$
 By continuity of $g_c(c,\cdot)$ and by continuity of trajectories  of $D^d$, we have $\tau_\varepsilon>0$.
 Then, by Proposition \ref{Prop:dynkin}\,(2) and taking into account the definition of $\tau_\varepsilon$, we have
 \beqs
 v_c(c,d)&=&\inf_{\sigma\in\mathcal{T}}\sup_{\tau\in\mathcal{T}} J(c,d;\sigma,\tau)\ \ \geq \ \ \inf_{\sigma\in\mathcal{T}} J(c,d;\sigma,\tau_\varepsilon)\\
 &=& \inf_{\sigma\in\mathcal{T}}\E\Big[\int_0^{\sigma\wedge\tau_\varepsilon} e^{-\rho t} g_c(c,D_t^d)\ud t+q_0^-e^{-\rho \sigma}\mathbf{1}_{\{\sigma<\tau_\varepsilon\}}-q_0^+e^{-\rho \tau_\varepsilon}\mathbf{1}_{\{\tau_\varepsilon<\sigma\}}\Big]
 \\
 &\geq& \inf_{\sigma\in\mathcal{T}}\E\Big[(\varepsilon-q_0^+)(1-e^{-\rho(\sigma\wedge\tau_\varepsilon)})+q_0^-e^{-\rho \sigma}\mathbf{1}_{\{\sigma<\tau_\varepsilon\}}-q_0^+e^{-\rho \tau_\varepsilon}\mathbf{1}_{\{\tau_\varepsilon<\sigma\}}\Big]\\
 &\geq& \inf_{\sigma\in\mathcal{T}} \E\Big[\varepsilon(1-e^{-\rho\tau_\varepsilon})\mathbf{1}_{\{\tau_\varepsilon<\sigma\}}-q_0^+e^{-\rho \tau_\varepsilon}\mathbf{1}_{\{\tau_\varepsilon<\sigma\}}\Big].
\enqs
Clearly the last term of the inequality above is larger than $-q_0^+$. Now, assume by contradiction that it is equal to $-q_0^+$. This means that there exists a minimizing sequence of stopping times $(\sigma_n)_{n\in\N}\subset \mathcal{T}$ such that
\beq\label{jjhq}
\lim_{n\rightarrow \infty} \, \E\Big[\varepsilon(1-e^{-\rho\tau_\varepsilon})\mathbf{1}_{\{\tau_\varepsilon<\sigma_n\}}-q_0^+e^{-\rho \tau_\varepsilon}\mathbf{1}_{\{\tau_\varepsilon<\sigma_n\}}\Big] \ \ =\ \  -q_0^+.
\enq
Hence, looking at the second addend in the expectation above, since the first one is nonnegative, we see that  we must have $\P\{\tau_\varepsilon<\sigma_n\}\rightarrow 1$. But then we must have 
$$(1-e^{\rho\tau_\varepsilon})\mathbf{1}_{\{\tau_\varepsilon<\sigma_n\}}\ \ \stackrel{\P}{\longrightarrow}\ \ 1-e^{-\rho \tau_\varepsilon}>0,$$ from which we deduce that 
$$\lim_{n\rightarrow\infty} \,\E\Big[\varepsilon(1-e^{-\rho\tau_\varepsilon})\mathbf{1}_{\{\tau_\varepsilon<\sigma_n\}}\Big]\ \ >\ \ 0,$$
contradicting \eqref{jjhq}. So we have shown that $v_c(c,d)>-q_0^+$. By continuity of $v_c(c,\cdot)$, this shows that $c>\hat{c}_+(d)$, completing the proof.\ep

\subsection{The dynamic programming equation}
The dynamic programming equation for the singular stochastic control problem \reff{tildev} takes the form of a  variational inequality:
\beq \label{VI00}
\max\big\{\, [\mathcal{L}v(c,\cdot)](d)-g(c,d),\ -v_c(c,d)-{q_0^+}, \ v_c(c,d)-q^-_0 \, \big\} &=&  0, \ \ \ (c,d)\in\,{{\mathcal{S}}},
\enq
where the second-order ordinary differential operator $\mathcal{L}$ is defined in \eqref{ODE}. 
Formally, \eqref{VI00} may be derived, assuming sufficient regularity of $v$ and exploiting its convexity in $c$, by looking at the three possibilities one has: (1) wait; (2) invest a small amount $\varepsilon$; (3) disinvest a small amount $\varepsilon$. We refer to \cite{FS}  for  a formal derivation of the dynamic programming equation in the  general context of singular control problems, and specifically to \cite{MZ} for a problem very similar to ours.

In the following, given a locally bounded function {$\phi:\mathcal{U}\rightarrow \mathbb{R}$}, where $\mathcal{U}\subset \R^n$ is an open set, we denote respectively by {$\phi^*$}, and  
{$\phi_*$}  the upper semicontinuous  and the lower semicontinuous envelope of $\phi$.  
  Since we do not know a priori if there exists a smooth solution to \reff{VI00}, we first rely in general on the notion of viscosity solutions:
 
 \begin{Definition}
 \noindent (i) We say that $v$ is a viscosity subsolution to \eqref{VI00} if  for any $(c,d)$ $\in$ $\Sc$, 
\beqs 
\max\big\{\, [\mathcal{L}\varphi(c,\cdot)](d)-g(c,d),\ -\varphi_c(c,d)-{q_0^+}, \ \varphi_c(c,d)-q^-_0\, \big\} & \leq &  0,
\enqs
whenever {$\varphi$ $\in$ $C^{1,2}(\mathcal{S}; \mathbb{R})$, $v^*(c,d)=\varphi(c,d)$,  and $v^*-\varphi$ has a local maximum at $(c,d)$.}

\noindent (ii)  We say that $v$ is a viscosity supersolution to \eqref{VI00} if  for any $(c,d)$ $\in$ $\Sc$, 
\beqs 
\max\big\{\,[\mathcal{L}\varphi(c,\cdot)](d)-g(c,d),\ -\varphi_c(c,d)-{q_0^+},\ \varphi_c(c,d)-q_0^-\,\big\} & \geq &  0,
\enqs
whenever {$\varphi$ $\in$ $C^{1,2}(\mathcal{S}; \mathbb{R})$, $v_*(c,d)=\varphi(c,d)$,  and $v_*-\varphi$ has a local minimum at $(c,d)$.}

\noindent  (iii) We say that $v$ is a viscosity solution to \eqref{VI00}  if {it} is both a viscosity sub- and supersolution.
\end{Definition}

\vspace{2mm}

The viscosity property of the value function follows usually from the dynamic programming principle (DPP).  The statement of DPP calls upon  delicate  measurable selection arguments.  Once we  know a priori that the value function is continuous,  {one can overcome this difficulty by exploiting the continuity}, see e.g. \cite{FS}. However, since the control set is unbounded,  and we are not assuming 
Lipschitz continuity of the coefficients in  \eqref{eq:demand} and - overall - of $g$, it is  not clear how  to get the continuity of the value function from its very definition. 
Instead, we can use the concept of weak dynamic programming introduced in \cite{BT}, which holds for our problem (see also Remarks 3.10 and 3.11 in \cite{BT}), stating that, for each $(c,d)\in\mathcal{S}$ and for each family $(\tau_I)_{I\in\mathcal{I}}$ of stopping times indexed by $I\in\mathcal{I}$, it holds
\begin{multline}\label{DPP}
 \inf_{I\in \mathcal{I}} \ \E\left[\int_0^{\tau_I^-}  e^{-\rho t} g(C_t^{c,I},D_t^d)\ud t+q_0^+\ud I_t^++q_0^-\ud I_t^- + e^{-\rho \tau_I}v_*(C_{\tau_I^-}^{c,I},D_{\tau_I}^d)\right]\\ \leq  \  v(c,d)\
 \leq \  \inf_{I\in \mathcal{I}} \ \E\left[\int_0^{\tau_I^-}  e^{-\rho t} g(C_t^{c,I},D_t^d)\ud t+q_0^+\ud I_t^++q_0^-\ud I_t^- + e^{-\rho \tau_I}v^*(C_{\tau_I^-}^{c,I},D_{\tau_I}^d)\right].  
\end{multline}

\begin{Proposition}
The value function $v$ is a viscosity solution to \eqref{VI00} on $\mathcal{S}$.
\end{Proposition}
\noindent
\textbf{Proof.} Given the weak DPP \eqref{DPP}, the proof is straightforward (and we omit it for brevity), and follows the line of the proof based on the standard Dynamic Programming Principle. Indeed, what one really needs are the two inequalities of \eqref{DPP} separately to prove the two viscosity properties separately. We can refer to \cite[Sec.\,5]{BT} where this is done for the case of continuous control; the proof can be adapted to our case of stochastic control.\ep


{\begin{Remark} 

{\rm  A comparison principle to the variational inequality \eqref{VI00} for viscosity sub-and super solution satisfying the growth condition \reff{growthv} could be proved using standard techniques (see \cite{CIL}), hence providing a uniqueness viscosity characterization of the value function $v$.  However, in our approach we rely mainly on  the viscosity property in order to derive a smooth-fit property. 
\ep  
}
\end{Remark}

\vspace{2mm}

We now investigate  the structure of the value function $v$ in the continuation region $\mathcal{C}$ and in the action regions $\mathcal{A}^\pm$.
The following lemma characterizes the structure of $v$ in the $c$-sections $S_c$ defined in \eqref{secc}.
\begin{Lemma}\label{prop:C2} Let  $c\in(\underline{c}_+ , \bar{c}_-)$.
\begin{enumerate}
 \item $v(c,\cdot)$ is a viscosity solution of the ODE
\beq \label{viscovCc}
[\mathcal{L}v(c,\cdot)](d)-g(c,d) &=& 0, \ \ \ \ d\in(\hat{d}_+(c),\hat{d}_-(c)).
\enq
\item
$
v(c,\cdot) \; \in  \;  C^2((\hat d_-(c),\hat{d}_+(c));\mathbb{R}).
$
\item  There exist constants $A(c),B(c)\in\mathbb{R}$ such that
\beq \label{vCc}
v(c,d) &=& A(c)\psi(d)+B(c)\varphi(d)+{{\hat{V}}}(c,d), \ \ \ \forall\,d\in (\hat{d}_-(c),\hat{d}(c)).
\enq
Moreover, \eqref{vCc} holds also at $\hat{d}_-(c), \hat{d}_+(c)$ when they do not coincide with $d_{\min}, d_{\max}$, respectively.
\end{enumerate}
\end{Lemma}
\noindent \textbf{Proof.} 1. Let us show the subsolution property (the proof of the supersolution property is completely analogous). 

First of all we note that, 
 since $v(\cdot,d)\in C^1(\R;\R)$, it is $v(c,d)=v(c_0,d)+\int_{c_0}^cv_c(\xi,d)\ud \xi$, for every $c,c_0\in\R$ and every $d\in\mathcal{O}$. Thus, since by Proposition \ref{prop:vc}\,(1) $v_c$ is  continuous in $\mathcal{S}$, we deduce the equalities
 \beq\label{unpp}
v^*(c,d)& =& v(c,\cdot)^*(d), \ \ \ \ \ \ \ \ \ \ \ \forall (c,d)\in\mathcal{S};\\
v^*(c,d)-v^*(c_0,d)  & =  & v(c,d)-v(c_0,d),  \ \ \ \ \forall (c,d)\in\mathcal{S}, \  \forall c_0\in\R.\label{op}
\enq

Let $c_0\in(\underline{c}_+ , \bar{c}_-)$, $d_0\in (\hat{d}_+(c_0),\hat{d}_-(c_0))$, and let $\phi\in C^2(\mathcal{O};\R)$ be such that 
\beq \label{visc}
\phi(d_0)&=& v(c_0,\cdot)^*(d_0), 
 \ \ \ \ \ \ \phi(d)\ \ \geq \  v(c,\cdot)^*(d), \ \ \forall d\in\mathcal{O}.
 \enq 
We claim that
\beq\label{super}
(v_c(c_0,d_0),\phi'(d_0), \phi''(d_0))&\in& D^{1,2,+}_{c,d} v^*(c_0,d_0),
\enq
where $D^{1,2,+}_{c,d}v^*(c_0,d_0)$ is the superdifferential of $v^*$ at $(c_0,d_0)$ of first order w.r.t. $c$ and of second order w.r.t. $d$ (see \cite{YZ}, Ch.\,4, Sec.\,5). We have to check that
\beq\label{ssuper}
\limsup_{(c,d)\rightarrow (c_0, d_0)} \!\!\! \frac{v^*(c,d)\!-\!v^*(c_0,d_0)\!-\! v_c(c_0,d_0)(c-c_0)\!
-\!\phi'(d_0)(d-d_0)\!-\!\phi''(d_0)(d-d_0)^2}{{|c-c_0|+|d-d_0|^2}} \ \leq \  0.
\enq
By \eqref{unpp} it has to be $(\phi'(d_0),\phi''(d_0))\in D^{2,+}_dv^*(c_0,d_0)$, where $D^{2,+}_{d}v^*(c_0,d_0)$ is the superdifferential of $v^*$ at $(c_0,d_0)$ of second order w.r.t. $d$.
Hence
\beq\label{ssuperd}
 {v^*(c_0,d)-v^*(c_0,d_0)
-\phi'(d_0)(d-d_0)-\phi''(d_0)(d-d_0)^2}& \leq & o(|d-d_0|^2).
\enq
Moreover, since $v(\cdot,d)\in C^1(\mathbb{R};\R)$ for every $d\in\Oc$ and $v_c$ is locally uniformly continuous w.r.t. $(c,d)\in\mathcal{S}$, for all  $\varepsilon>0$ there exists $\delta>0$ such that 
\beq
{v(c,d)-v(c_0,d)
-v_c(c_0,d)(c-c_0)}\!\!& \leq &\!\! o(|c-c_0|), \, \mbox{{unif.\,\,in}}
\,  d\in(d_0-\delta,d_0+\delta),\label{rrf}\\
|v_c(c_0,d)-v_c(c_0,d_0)|&\leq &\varepsilon, \ \ \ \ \ \ \ \ \ \ \ \ \ \ \  \ \forall \,  d\in(d_0-\delta,d_0+\delta).\label{rff1}
\enq
By \eqref{op}, we  derive from \eqref{rrf}
\beq\label{rff2}
 \!\!{v^*(c,d)-v^*(c_0,d)
-v_c(c_0,d)(c-c_0)} &\leq & o(|c-c_0|),  \ \mbox{{unif.\,\,in}}
\,  d\in(d_0-\delta,d_0+\delta).
\enq
By subtracting and adding $v_c(c_0,d_0)(c-c_0)$ in \eqref{rff2} and using \eqref{rff1}, we get
\beq\label{ssuperc}
\!\!\!{v^*(c,d)\!-\!v^*(c_0,d)\!
-\!v_c(c_0,d_0)(c-c_0)} \leq  o(|c-c_0|)\!+ \!\varepsilon\! \cdot\!|c-c_0|, \, \mbox{{unif.\,\,in}}
\,  d\in(d_0-\delta,d_0+\delta).
\enq
Combining \eqref{ssuperd} and \eqref{ssuperc}, dividing by $|c-c_0|+|d-d_0|^2$, and taking the limsup, since $\varepsilon$ was arbitrary, we finally get \eqref{ssuper}, thus \eqref{super}. 

Now, starting from \eqref{super}, we can construct  (see, e.g.,  \cite{YZ}, Ch.\,4, Lemma 5.4\footnote{The proof works even if the function is just upper semicontinuous.}) a function $\varphi\in C^{1,2}(\mathcal{S};\R)$ such that $\varphi(c_0,d_0)=v^*(c_0,d_0)$, $\varphi\geq v^*$ on $\mathcal{S}$ and 
\beq\label{eqq}
(\varphi_c(c_0,d_0),\varphi_d(c_0,d_0), \varphi_{dd}(c_0,d_0))&=& (v_c(c_0,d_0),\phi'(d_0), \phi''(d_0)).
\enq
Now
notice that $-q_0^+< v_c(c_0,d_0)<q_0^-$, as $(c_0,d_0)\in\,\,{\mathcal{C}}$ (Proposition \ref{sezioni}). Hence, since $v$ is a viscosity solution to \eqref{VI00}, taking into account \eqref{eqq} we finally get the desired inequality $[\mathcal{L}\phi](d_0)\leq 0$.
\smallskip

2. Let $c\in(\underline{c}_+ , \bar{c}_-)$ and, given $a,b\in \bar{S_c}$ with $a<b$, consider the Dirichlet problem
\begin{equation}\label{diffpr}
\begin{cases}
\rho u(d)-\mu(d)u'(d)-\frac{1}{2}\sigma^2(d)u''(d)\ =\ g(c,d), \ \ \ \ d\in(a,b),\\
u(a)\ =\ v(c,a), \ \ \ u(b)\ =\ v(c,b).
\end{cases}
\end{equation}
{This problem clearly admits a unique viscosity solution, which must coincide with $v(c,\cdot)$ in $[a,b]$ by item 1. 
On the other hand, by  since $\sigma^2(\cdot)>0$, \eqref{diffpr} is a uniformly elliptic  problem, so it   admits a solution of class $C^0([a,b];\mathbb{R})$ $\cap$ 
$C^2((a,b);\mathbb{R})$, which is also a viscosity solution, so coincides with $v$.} Hence, we deduce that  $v(c,\cdot)\in  C^2((\hat{d}_-(c),\hat{d}(c));\mathbb{R})$, and 
satisfies in a classical sense:
\beqs 
[\mathcal{L} v(c,\cdot)](d) - g(c,d) &=& 0,  \;\;\;\;\;   d \in (\hat d_-(c),\hat d_+(c)). 
\enqs

3. Notice that  $\hat V(c,\cdot)$ is a particular solution to the ODE 
\beq \label{nonhomODE}
[\mathcal{L} \phi(c,\cdot)](d) - g(c,d) &=& 0, \ \ \ \ d\in \Oc.
\enq
Therefore the general solution to \reff{nonhomODE} is  in the form:
\beqs
A(c)\psi(d)+B(c)\varphi(d)+{{\hat{V}}}(c,d), \ \ \ \ d\in S_c
\enqs
for some real-valued constants $A(c)$, $B(c)$, which  proves, together with item 2,  the structure   
\reff{vCc} of $v$ in $S_c$.

The extension of  \eqref{vCc} at $\hat{d}_-(c)$ and at $\hat{d}(c)$, when they do not coincide with $d_{\min}, d_{\max}$, respectively, can be obtained by taking $a=\hat{d}_-(c)$ and $b=\hat{d}_+(c)$ in the argument above.
\ep

\vspace{2mm}

\begin{Lemma}\label{meninf} 
We have
\beq\label{plc}
\lim_{d\downarrow d_{\min}} (v(c,d)-{\hat{V}}(c,d))& =& 0,\ \  \forall c\in(\underline{c}_+,\underline{c}_-);\\
\lim_{d\uparrow d_{\max}} (v(c,d)-{\hat{V}}(c,d))&=&0,\ \  \forall c\in(\bar{c}_+,\bar{c}_-).\label{plc1}
\enq
\end{Lemma}
\noindent \textbf{Proof.}   We prove \eqref{plc}, the proof of \eqref{plc1} is  analogous.

Fix $c\in(\underline{c}_+,\underline{c}_-)$. In this case we have $\hat{d}_-(c)=d_{\min}$. Then, due to Lemma  \ref{prop:C2}, we have that  $v(c,\cdot)$ $\in$ $C^2(({d}_{min},\hat{d}_+(c));\mathbb{R})$, and that it satisfies in a classical sense
\beq \label{pfe}
[\mathcal{L}v(c,\cdot)](d) -  g(c,d) &=& 0, \ \ \ \ \forall d\in (d_{\min},\hat{d}_+(c)).
\enq
Let 
 $d_0\in(d_{\min},\hat{d}_+(c))$ be fixed and take a generic $d\in(d_{\min},d_0)$. Consider the stopping time 
 $$\tau_d \ \ = \ \ \inf\,\{t\geq 0 \ | \ D_t^d\geq d_0\}.$$
Since $d_{\min}$ is  not-entrance for the diffusion $D$, we have  
(see e.g. \cite[Ch.\,20]{Kall}):
\beq \label{qws}
\tau_d \nearrow \infty   &   \mbox{ when } \;\;\; d\downarrow d_{\min}.
\enq
Given a sequence $(d_n)\subset (d_{\min},d)$ such that $d_n\downarrow d_{\min}$ consider the stopping times 
 $$\tau_d^n \ \ = \ \ \inf\,\{t\geq 0 \ | \ D_t^d\leq d_n\}.$$
Since  $d_{\min}$ is inaccessible for the diffusion $D$, we have
\beq \label{qws1}
\tau_d^{n} \nearrow \infty   &    \mbox{ when } \;\;\; n\rightarrow \infty.
\enq
By \eqref{pfe} and  definition of $\tau_d$, we apply It\^o's formula  to $v(c,D_t^d)$ in the interval $[0,\tau_d\wedge \tau_d^n\wedge n)$, 
\beqs
v(c,d) &=&\int_0^{\tau_d\wedge \tau_d^n\wedge n}e^{-\rho t} g(c,D^d_t)\ud t+\int_0^{\tau_d\wedge \tau_d^n \wedge n}e^{-\rho t} v_d(c,D^d_t)\ud W_t+ e^{-\rho \tau_d} v(c,D_{\tau_d\wedge \tau_d^n \wedge n}^d).
\enqs
By taking the expectation (noting that the expectation of the stochastic integral vanishes by our localization and that $v\geq 0$), we get 
\beqs
v(c,d) &\geq & 
\mathbb{E}\Big[\int_0^{\tau_d\wedge \tau_d^n\wedge n}e^{-\rho t} g(c,D^d_t)\ud t\Big].
\enqs
By taking the limit for $n \rightarrow\infty$ (note that $g\geq 0$, so we can use monotone convergence) and using \eqref{qws1}, we get
\beqs
v(c,d) &\geq& \mathbb{E}\Big[\int_0^{\tau_d}e^{-\rho t} g(c,D^d_t)\ud t\Big].
\enqs
Subtracting ${\hat{V}}(c,d)$ in both sides of the inequality above, we get
\beqs
v(c,d)-\hat{V}(c,d) & \geq &  \mathbb{E}\Big[\int_{\tau_d}^\infty e^{-\rho t} g(c,D^d_t)\ud t\Big]
\enqs
Taking the liminf for $d\downarrow d_{\min}$, and using \eqref{qws}, we obtain
\beqs
\liminf_{d\downarrow d_{\min}} \ (v(c,d)-\hat{V}(c,d)) &\geq& 0,
\enqs
and so the required limiting result, since we always have $v$ $\leq$ $\hat V$ (see \reff{mingamma}).\ep

%
%
%
%
%

 \subsection{Structure of the value function} 
 We  can now provide the complete structure of the value function. Let us define 
 \beqs
\Oc_+&:=& \{d\in\mathcal{O} \ | \ {\hat{c}_+}(d)>-\infty\}, \ \ \ \ \Oc_-\ \ :=\ \  \{d\in\mathcal{O} \ | \ {\hat{c}_-}(d)<\infty\}.
 \enqs
  Note that $\Oc_\pm$ are connected due to monotonicity of $\hat{c}_\pm$.

 \begin{Theorem}  \label{prop:structure} 
 (Structure and properties of the value function) 
  
\noindent There exist functions  $$A,B \in\  C^1((\underline{c}_+,\bar{c}_-);\R), \ \ \ \ \ z_{\pm}:  
\Oc_{\pm} \rightarrow \mathbb{R},$$
(with $A,B$ eventually extendable to $C^1$ functions up to $\underline{c}_+,\bar{c}_-$, respectively, when there exists $d\in\mathcal{O}$ such that ${\hat{c}_+}(d)=\underline {c}_+$, or when there exists $d\in\mathcal{O}$ such that ${\hat{c}_-}(d)=\bar {c}_-$),  such that  
\begin{equation}\label{VII1}
 v(c,d) \; = \; \begin{cases}A(c)\psi(d)+B(c)\varphi (d)+{{\hat{V}}}(c,d), \ \ \ \ \ \  \mbox{on} \ \bar{\mathcal{C}},\\
z_+(d) -{q_0^+}c, \ \ \ \ \ \ \ \ \ \ \ \  \ \ \ \ \ \ \ \ \ \ \ \ \ \ \ \ \, \mbox{on} \ \mathcal{A}^+,\\
 z_-(d) +{q_0^-}c, \ \ \ \ \ \ \ \ \ \ \ \  \ \ \ \ \ \ \ \ \ \ \ \   \, \ \ \ \ \mbox{on} \ \mathcal{A}^-.
 \end{cases}
 \end{equation}
Moreover:
\begin{itemize}
\item[(i)] 
$A(c)\ = \ 0$ for every $c\in[\bar{c}_+,\bar{c}_-)$, and $B(c)\ = \ 0$ for every  $c\in (\underline{c}_+,\underline{c}_-]$ (note that these intervals may be empty).

\item[(ii)] $z_{\pm}$  can be written in terms of the values of $v$ at $\partial \mathcal{C}$ and of ${\hat{c}_{\pm}}$ as
\beq\label{expz}
z_+(d) \;=\;  v({\hat{c}_+}(d),d)+{q_0^+}{\hat{c}_+}(d), \ \ \ \ d\in\Oc_+,\\
 z_-(d)  \;=\;  v({\hat{c}_-}(d),d)-{q_0^-}{\hat{c}_-}(d), \ \ \  \ d\in\Oc_-.
\enq
\end{itemize}
\end{Theorem}
\noindent \textbf{Proof.}
\emph{Structure of $v$ in $\bar{\mathcal{C}}$.} By Lemma \ref{prop:C2}(3), we already know that there exist functions $A,B: (\underline{c}_+,\overline{c}_-)\rightarrow \R$ such that we have
\beq\label{structure}
v(c,d) &=& A(c)\psi(d)+B(c)\varphi(d)+\hat{V}(c,d), \ \  \ \ (c,d)\in\Cc.
\enq
Let $c_0\in (\underline{c}_+,\overline{c}_-)$. Since $\mathcal{C}$ is open, from the representation \eqref{cupcirc} we see that we can find $d,d_0\in\Oc$ such that $(c,d_0),(c,d)\in S_c$ for every $c\in(c_0-\varepsilon,c_0+\varepsilon)$, for some $\varepsilon>0$.
Writing    \eqref{structure} at $(c,d), (c,d_0) \in \mathcal{C}$, and taking into account that $\psi(d)\varphi(d_0)-\varphi(d)\psi(d_0)\neq 0$ for all $d\neq d_0$ (this is due to strict monotonicity of $\varphi, \psi$), we can retrieve $A,B$ in the interval $(c_0-\varepsilon,c_0+\varepsilon)$ as 
\beq\label{AB}
A(c)&=&\frac{(v(c,d)-\hat{V}(c,d))\varphi(d_0)-(v(c,d_0)-\hat{V}(c,d_0))\varphi(d)}{\psi(d)\varphi(d_0)-\varphi(d)\psi(d_0)},\\
B(c)&=&\frac{(v(c,d_0)-\hat{V}(c,d_0))\psi(d)-(v(c,d)-\hat{V}(c,d))\psi(d)}{\psi(d)\varphi(d_0)-\varphi(d)\psi(d_0)}.\label{AB2}
\enq
Hence, since $v(\cdot,d)$ and $\hat{V}(\cdot,d)$ are of class $C^1$ for any fixed $d\in\mathcal{O}$,  we get, by arbitrariness of $c_0$, that  $A,B\in C^1((\underline{c}_+,\bar{c}_-);\R)$.

Now assume that there exists $d\in\mathcal{O}$ such that ${\hat{c}_+}(d)=\underline {c}_+$. Then, since the function $\hat{c}_+$ is nondecreasing and right-continuous, there exists an interval $(a,b)\subset\Oc$ such that ${\hat{c}_+}(d)=\underline {c}_+$ in $(a,b)$. Take $d_0,d\in(a,b)$. Then, for every   $c>\underline{c}_+$, it is $(c,d_0), (c,d)\in\Cc$. We can then write the relation  \eqref{AB} for every $c> \underline{c}_+$ and pass it  to the limit for $c\downarrow \underline{c}_+$. In such a way we see that  $A$ can be extended to $C^1$ function up to $\underline{c}_+$. The same argument holds true for the other case involving $B$ and $\bar{c}_-$.
\vspace{1mm}

Let us  now check that   \eqref{structure}  also holds at the points of the boundary $\partial\mathcal{C}$.
Let $(c,d)\in\partial^+ \mathcal{C}$. In this case,  one of the following case must hold\,: 
\begin{itemize}
\item[(a)] $d=\hat{d}_+(c)\in\Oc$, 
\vspace{-.15cm}
\item[(b)] $c={\hat{c}_+}(d)$ and $\{(c,d) \ | \ c\in({\hat{c}_+}(d), \hat{c}(d)+\varepsilon)\}\subset \ {\mathcal{C}}$ for some $\varepsilon>0$,
\vspace{-1mm}
\item[(c)] $d=\hat{d}_+(c')$ for $c'\in (c,c+\varepsilon)$ for some $\varepsilon>0$.
\end{itemize}
In the case (a) the form  \eqref{structure} holds by Lemma \ref{prop:C2}\,(3). In the case (b) the structure \eqref{structure} holds  by continuity of  $A,B$ and of $v$ with respect to $c$, and by the already proved structure in  ${\mathcal{C}}$. In the case (c) the structure  \eqref{structure} holds by case (a) and by continuity of $A, B$ and of  $v$ with respect to $c$. 

The same argument holds for points belonging to the boundary $\partial^-\mathcal{C}$, so we conclude that
\beq \label{structureclosure}
v(c,d) &=& A(c)\psi(d)+B(c)\varphi(d)+{{\hat{V}}}(c,d), \ \ \ \ \mbox{in} \ \bar{\mathcal{C}}.
\enq

\smallskip  

\noindent \emph{Structure of $v$ in $\mathcal{A}^{\pm}$.} This follows directly from the  definition \eqref{mathD} of $\mathcal{A}^{\pm}$. 

\medskip

Let us now prove the remaining properties. 

\noindent \emph{(i)} 
 Let $c\in(\bar{c}_+,\bar{c}_-)$. We can use \eqref{structureclosure} and write  
\beqs
\lim_{d\uparrow d_{\max}} v(c,d)&=& \lim_{d\uparrow d_{\max}} (A(c)\psi(d)+B(c)\varphi (d)+\hat{V}(c,d)).
\enqs
By taking into account  Lemma \ref{meninf} and \reff{psiphi},  we see that it must be  $A(c)$ $=$ $0$. In a similar way one proves that $B(c)=0$ for every $c\in(\underline{c}_+,\underline{c}_-)$.  Then $A(\bar{c}_+)=0$ and $B(\underline{c}_-)=0$ follow by continuity.

\medskip

\noindent \emph{(ii)} It  follows using \eqref{VII1} and by evaluating $v$  at the points $({\hat{c}_{\pm}}(d),d)\in\bar{\mathcal{C}}$.
%
%
\ep
\subsection{Optimal control}

 In the following we suppress, for simplicity of notation, the superscript $d$ in $D^d$. Moreover, the superscript $k$ in the notation  $C^k_t$  below will not denote the initial datum, but a running natural index.

Let $(c,d)\in\mathcal{S}$. Let us define,
with the convention $\inf\emptyset=\infty$, the random times 
\beqs
{\tau}_0^+& := & \inf \ \{ t\geq 0 \ | \ c<\hat{c}_+(D_t)\}, \  \ \ \ {\tau}_0^-\ \ := \ \ \inf \ \{ t\geq 0 \ | \ c>\hat{c}_-(D_t)\}, \ \ \ \ \tau_0\ \ :=\ \ \tau_0^+\wedge\tau_0^-.
\enqs
Due to \eqref{c+-}, we have $\{\tau_0^+=\tau_0^-\}=\{\tau_0=\infty\}.$ Define also
$$\Omega_\infty\ \ :=\ \ \{\tau_0=\infty\}, \ \ \ \ \ \Omega_+\ \ : =  \ \{\tau_0^+<\tau_0^-\}, \ \ \ \ \ \ \Omega_-\ \ : =  \ \{\tau_0^+>\tau_0^-\}.$$
Define
\beqs
C_t^{0}\ \ =\ \ c, \ \ \ \ \ t\geq 0,
\enqs
and define recursively the following processes and stopping times\,:
\begin{itemize}
\item[-] For all $k\geq 0$,
\beqs
\overline{D}^k_t& := & \max_{s\in[\tau_{k-1},t]} D_s, \ \ \ \ \ \ \underline{D}^k_t\ \ := \ \  \min_{s\in[\tau_{k-1},t]} D_s, \ \ \ \ \ \ t\geq \tau_{k-1},
\enqs
\item[-] If $k\geq 1$ is odd, 
\beqs
C_t^{k}\ \ :=\ \ \begin{cases}
c, \ \ \ \ \ \ \ \ \ \ \ \ \ \ \ \ \ \ \mbox{on} \ \Omega_\infty,\\
c+\hat{c}_+(\overline{D}^k_t), \ \ \ \ \ \mbox{on} \ \Omega_+,\\
c+\hat{c}_-(\underline{D}^k_t), \ \ \ \ \ \mbox{on} \ \Omega_-,\\
\end{cases}
\ \ \ \ \ \ \ t\ \geq \ \tau_{k-1},
\enqs
\beqs
\tau_k\ \ :=\ \ \begin{cases}
\infty, \ \ \ \ \ \ \ \ \ \ \ \ \ \ \ \ \ \ \ \ \ \ \ \ \ \ \ \ \ \ \ \ \ \ \ \ \ \ \ \  \ \ \mbox{on} \ \Omega_\infty,\\
\inf\ \{t\geq \tau_{k-1} \ | \ C^{*,k}_t>\hat{c}_-({D}_t)\}, \ \ \ \ \ \mbox{on} \ \Omega_+,\\
\inf\ \{t\geq \tau_{k-1} \ | \ C^{*,k}_t<\hat{c}_+({D}_t)\}, \ \ \ \ \ \mbox{on} \ \Omega_-.
\end{cases}
\enqs
\item[-] If $k\geq 2$ is even
\beqs
C_t^{k}\ \ :=\ \ \begin{cases}
c, \ \ \ \ \ \ \ \ \ \ \ \ \ \ \ \ \ \ \mbox{on} \ \Omega_\infty,\\
c+\hat{c}_+(\overline{D}^k_t), \ \ \ \ \ \mbox{on} \ \Omega_-,\\
c+\hat{c}_-(\underline{D}^k_t), \ \ \ \ \ \mbox{on} \ \Omega_+,\\
\end{cases}
\ \ \ \ \ \ \ t\ \geq \ \tau_{k-1},
\enqs
\beqs
\tau_k\ \ :=\ \ \begin{cases}
\infty, \ \ \ \ \ \ \ \ \ \ \ \ \ \ \ \ \ \ \ \ \ \ \ \ \ \ \ \ \ \ \ \ \ \ \ \ \ \ \ \  \ \ \mbox{on} \ \Omega_\infty,\\
\inf\ \{t\geq \tau_{k-1} \ | \ C^{*,k}_t>\hat{c}_-({D}_t)\}, \ \ \ \ \ \mbox{on} \ \Omega_-,\\
\inf\ \{t\geq \tau_{k-1} \ | \ C^{*,k}_t<\hat{c}_+({D}_t)\}, \ \ \ \ \ \mbox{on} \ \Omega_+.
\end{cases}
\enqs
\end{itemize}
%
%
Since  $\mathcal{A}^\pm$ are closed and  $\sigma^2>0$, we have, if $k$ is odd
\beqs
\inf\,\{t\geq \tau_k \ | \ (C^{*,k}_t,D_t)\in \,\stackrel{\circ}{\Ac^+}\} & =& \inf\,\{t\geq \tau_k \ | \ C^{*,k}_t<\hat{c}_+(D_t)\}, \ \ \mbox{a.e. in } \Omega_+,\\
\inf\,\{t\geq \tau_k \ | \ (C^{*,k}_t,D_t)\in \, \stackrel{\circ}{\Ac^-}\} & =& \inf\,\{t\geq \tau_k \ | \ C^{*,k}_t>\hat{c}_-(D_t)\}, \ \ \mbox{a.e. in } \Omega_-,
\enqs
and similar representations if $k$ is even. 
Hence, since $\F$ satisfies the usual conditions, so hitting times of open sets are stopping times, we see that the sequence $(\tau_k)$ is a sequence of stopping times. 
%
%
\vspace{2mm}
Setting $\tau_{-1}:=0$, define the process
\beq\label{optstate}
C^*_t&:=& \sum_{k=0}^\infty C^{k}_t\, \mathbf{1}_{[\tau_{k-1}, \tau_k)}(t), \ \ \ \ \ \ \ \ t\geq 0. 
\enq
{Since  $\tau_k\rightarrow \infty$ almost surely, the process $C^*$ is  well defined for every $t\geq 0$. Moreover it is clearly right-continuous and adapted. By construction
\beq\label{CinC}
(C^*_t,D_t)&\in&\bar{\Cc}, \ \ \ \ \forall t\geq 0.
\enq
Define the control 
\beq\label{optI}
I^*_t\ \ :=\ \ C^*_t-c.
\enq
The control process $I^*$  does the minimum effort to keep the couple $(C_t^{*},D_t)$ inside $\bar{\mathcal{C}}$. 
More precisely, at time $t\geq 0$:
\begin{itemize}
\item[-] if $(C_{t^-}^*,D_t)\in\Cc$, no action is taken ($\ud I^*=0$);
\item[-] if  $(C_{t^-}^*,D_t)\in\partial\Cc$ (e.g., assume $(C_{t^-}^*,D_t)\in\partial^+ \Cc$; simmetrically one can argue in the case $(C_{t^-}^*,D_t)\in\partial^- \Cc$), then two cases have to be distinguished:
\begin{itemize}
\item[-] if $C_{t^-}^*=\hat{c}_+(D_t)$ (which occurs in particular if $\hat{c}$ is continuous at $D_t$), then $I^*$ acts in order to reflect  $(C^*_t,D_t)$  at the boundary $\partial \Cc^+$  along the positive $c$-direction. Note that no action is taken if $\hat{c}_+$ is constant in a right-neighborhood of $D_t$.
\item[-] if $\hat{c}_+$ is discontinuous at $D_t$ and $C_{t^-}^*<\hat{c}_+(D_t)$, then the process $C^*$ has a positive jump  $ \Delta C^*_t= \Delta I^{*,+}_t= \hat{c}_+(D_t)- C^*_{t^-}$. 
\end{itemize}
\end{itemize}
Regarding the last possibility, letting $\mathcal{N}^{\pm}$ be the (at most countable) sets of discontinuity points of ${\hat{c}_{\pm}}$, respectively, due to the continuity of trajectories of ${D}$, we see that the process $I^*=I^{*,+}-I^{*,-}$ can jump
\begin{itemize}
\item[(a.1)] either at time 
$0$ when $c<{\hat{c}_+}(d)$ or when $c>\hat{c}_-(d)$, and in this case we have, respectively, $\Delta I^*_0=\Delta I_0^{*,+}= {\hat{c}_+}(d)-c$ or $\Delta I^*_0=-\Delta I_0^{*,-}= {\hat{c}_-}(d)-c$;
\item[(a.2)] when $D_t\in\mathcal{N}^+$ and $C^{*}_{t^-}<{\hat{c}_+}(D_t)$, and  in this case $\Delta I^*_t=\Delta I^{*,+}_t ={\hat{c}_+}(D_t)-C^{*}_{t^-}$. 
\item[(a.3)] when $D_t\in\mathcal{N}^-$ and $C^{*}_{t^-}>{\hat{c}_-}(D_t)$, and  in this case $\Delta I^*_t=-\Delta I_t^{*,-}=C^{*}_{t^-} - {\hat{c}_-}(D_t)$. 
\end{itemize}}
\begin{Lemma}
The processes $C^*,I^*$ satisfy
\beq
\label{int}
 \int_0^{\infty}e^{-\rho t} \textbf{1}_{\{(C^{*}_t,D_t)\, \in \, \mathcal{C}\}}\,\ud I^{*,\pm}_t & =&  0. 
\enq
\end{Lemma}
\noindent \textbf{Proof.}
Fix $\omega\in \Omega$ and suppose that $(C^{*}_t(\omega),D_t^d(\omega))\in\mathcal{C}$. Then, by definition of the $\tau_k$'s and since $\Cc$ is open, we must have $t\in (\tau_{k-1}(\omega), \tau_{k}(\omega))$ for some $k\geq 0$, and 
\beq
C^*_t(\omega)\ \ \in\ \ \big(\,{\hat{c}_+}(D_t(\omega)), \ \hat{c}_-(D_t(\omega))\,\big).
\enq
By definition of $C^*$, $\tau_{k-1}, \tau_k$, we see that $C.^*(\omega)$ is constant in some suitable  neighborhood $(t-\varepsilon(\omega), t+\varepsilon(\omega))$ of $t$, hence also $I.^*(\omega)$ is constant therein. Thus, we have proved \eqref{int}.
\ep
\vspace{2mm}

The second main  result provides the existence and  an explicit  description of the optimal state process (and a description of the optimal investment in terms of the optimal state).

\begin{Theorem}\label{Verification}
(Optimal control)
\noindent  Let $(c,d)\in {\mathcal{S}}$.  The process $C^*$ constructed before  in \eqref{optstate} is 
an optimal state process for the value function at  $(c,d)$, with corresponding optimal control $I^*=(I^{*,+},I^{*,-})$ defined by \eqref{optI}.
\end{Theorem}
\noindent \textbf{Proof.} 
Let us show  that 
\beq \label{ine1}
v(c,d) &\geq&   \mathbb{E}\Big[\int_0^{\infty}e^{-\rho t} \Big(g(C_t^*,D_t)+{q_0^+}\ud I^{*,+}_t-q_0^-\ud I_t^{*,-}\Big)\Big]. 
\enq
Let $(K_n)$ be  an increasing  sequence of compact subsets of $\mathcal{S}$ such that $\cup_{n\in\mathbb{N}}K_n=\mathcal{S}$.
Consider the (bounded) stopping time $\tau_n$ $=$ $\inf\{t\geq 0 \ | \ C_t^*\wedge D_t\notin K_n \} \wedge n$,  and notice that 
$\tau_n\nearrow \infty$ a.s. when $n$ goes to infinity.   From \eqref{structureclosure} and since $\hat{V}\in C^{1,2}(\mathcal{S};\mathbb{R})$, we see that  $v\in C^{1,2}(\bar{\mathcal{C}};\mathbb{R})$. Thus, by  \eqref{CinC}, we may apply   
 It\^o's formula (see Proposition \ref{prop:ITO}) to $e^{-\rho t} v(C_t^*,D_t^d)$  between $0$ and $\tau_n$, take expectation, and obtain (after observing  that  the stochastic integral over the  interval $[0,\tau_n\wedge T)$ vanishes in expectation {due to our localization}): 
\beq
v(c,d) &=& \mathbb{E}\Big[e^{-\rho\tau_n}v(C^*_{\tau_n\wedge T},D_{\tau_n\wedge T})\Big] 
\;   +  \;  \mathbb{E}\Big[\int_0^{\tau_n}e^{-\rho t}[\mathcal{L}v(C_t^*,\cdot)](D_t){{\ud}t}\Big] \label{Ito} \\
& &  -   \; \mathbb{E}\Big[\int_0^{\tau_n}e^{-\rho t} v_c(C_t^*,D_t)\ud I^{*}_t\Big] \;  \nonumber \\
& & - \; \mathbb{E}\Big[\sum_{0\leq t\leq \tau_n}e^{-\rho t}(v(C_t^{*},D_t)-v(C_{t^-}^{*},D_t)-v_c(C_t^*,D_t) \Delta C_t^*)\Big], \nonumber 
\enq
Now observe that 
  $[\mathcal{L}v(c',\cdot)](d')=g(c',d')$ for $(c',d')$ in ${\mathcal{C}}$ but also in $\bar{\mathcal{C}}$ by continuity of $g$ and 
  {since $v\in C^{1,2}(\bar{\Cc};\R)$}.  This implies 
  \beq\label{hhg}
  \mathbb{E}\Big[\int_0^{\tau_n}e^{-\rho t}[\mathcal{L}v(C_t^*,\cdot)](D_t^d){{\ud}t}\Big]&=&  \mathbb{E}\Big[\int_0^{\tau_n}e^{-\rho t} \big(g(C_t^*,D_t^d)\ud t\Big].
  \enq
Now,   notice that $\ud I^{*,+} = 0$ if $(C_t^*,D_t^d)\in\mathcal{A}^-$ and $\ud I^{*,-} = 0$ if $(C_t^*,D_t^d)\in\mathcal{A}^+$. Then  
  taking into account \eqref{int} and the fact that $v_c=-{q_0^+}$ in $\mathcal{A}^+$ and $v_c=q_0^-$ in $\Ac^-$,  we have
\beq\label{hhg1}
-\ \mathbb{E}\Big[\int_0^{\tau_n}e^{-\rho t} v_c(C_t^*,D_t^d)\ud I^{*}_t\Big] &=&  \; \mathbb{E}\Big[ \int_0^{\tau_n} e^{-\rho t} (q_0^+ \ud I^{*,+}_t+q_0^- \ud I^{*,-}_t)\Big].
\enq
Moreover, considering the three possibilities of jump (a.1)--(a.3) described above for $I^*$, we have 
\beq\label{hhg2}
v(C_t^{*},D_t^d)-v(C_{t^-}^{*},D_t^d)-v_c(C_t^*,D_t^d) \Delta C_t^*&=&0, \ \ \ \ \forall t\geq 0.
\enq
Therefore by nonnegativity of $v$ and  \eqref{Ito}-\,-\eqref{hhg2}, we have
\beqs
v(c,d) &\geq& \mathbb{E}\Big[\int_0^{\tau_n}e^{-\rho t} \big(g(C_t^*,D_t^d)\ud t+{q_0^+}\ud I^{*,+}_t+q_0^-\ud I_t^{*,-}\big)\Big].
\enqs
Letting $n\rightarrow\infty$, from monotone convergence we get the inequality \eqref{ine1}. Since the opposite inequality always holds by definition of $v$, this proves the equality, i.e. that $I^*$ is an optimal control.
\ep

\vspace{2mm}

The picture below  represents a possible shape of the solution. The state space region $\mathcal{S}$ is the half-plane  on the right of the vertical dotted line. When the system lies in the continuation region $\mathcal{C}$, it moves along the horizontal lines and no action is taken. Whenever the system touches the boundary $\partial\mathcal{C}$, the optimal control (acting along the vertical lines as indicated by the arrows in the picture) consists in doing the minimal effort to keep the system in $\bar{\mathcal{C}}$. We notice that, if the boundary $\hat{c}_+$ or the boundary $\hat{c}_-$ is constant somewhere, no action is taken if the system reaches this part of boundary, and the system lies on this part of the boundary for a certain time until it meets a strictly increasing part of this boundary.
\begin{Remark}\label{Rem:positive}
{\rm
From the solution found, it turns out that when $\underline{c}_-\geq 0$, starting from $c\geq 0$ the optimal state process verifies $C^*\geq 0$. This means that the solution is, henceforth, also thee solution of the problem with state constraint $C\geq 0$.}
\end{Remark}

\begin{Corollary}\label{cor:fin}
\begin{enumerate}
\item 
If $\lim_{c\downarrow-\infty}g_c(c,d)=-\infty$, then $\hat{c}_+>-\infty$ in $(d,d_{\max})$. 
\item If $\lim_{c\uparrow\infty}g_c(c,d)=\infty$, then $\hat{c}_-<\infty$ in $(d_{\min},d)$. 
\end{enumerate}
\end{Corollary}
\noindent \textbf{Proof.} We prove item 1, then item 2 can be proved symmetrically.

Let $d\in\Oc$ be such that $\lim_{c\downarrow-\infty}g_c(c,d)=-\infty$. Take $c_0\in\R$ such that $g_c(c_0,d)\leq 0$ and $\hat{c}_-(d)> c_0$. Since by Assumption \ref{ass:cost}  $g_c$ is nondecreasing in $c$ and nonincreasing in $d$, we have $g_c\leq 0$ in $(-\infty,c_0]\times [d,d_{\max})$. Assume, by contradiction, that there exists $d_1\in(d,d_{\max})$ such that $\hat{c}_+(d_1)=-\infty$. By monotonicity of $\hat{c}_+$ this implies that $\hat{c}_+\equiv -\infty$ in $(d_{\min},d_1]$. 
Now, given any $c\leq c_0$ and $d_0\in(d,d_1)$, define the stopping times 
$$
\sigma =\inf\,\{t\geq 0 \ | \ D^{d_0}_t\leq d\}, \ \ \ \tau=\inf\,\{t\geq 0 \ | \ D_t^{d_0}\geq d_1\}, \ \  \ \tau^*(c)=\inf\,\{t\geq 0 \ | \ D^{d_0}_t\geq \hat{d}_+(c)\}.
$$
 Observe that $\tau\leq \tau^*(c)$, for every $c\in\R$, since $\hat{d}_+(c)$ has to be larger than $d_1$, as $\hat{c}_+\equiv-\infty$ in $(d_{\min},d_1]$. Moreover, by Proposition \ref{Prop:dynkin} and Theorem \ref{Verification}, $\tau^*(c)$ is the optimal stopping time of P2 for the Dynkin game defined in Subsection \ref{sub:dynkin}.
Hence, we must have, taking also into account that $g_c(c,\cdot)$ is nonincreasing, that  $g_c\leq 0$ in $(-\infty,c_0]\times [d,d_{\max})$, and that $\tau\leq \tau^*(c)$, 
\beqs
v_c(c,d) &\leq &J(c,d;\sigma,\tau^*(c))\\
&=& \E\left[\int_0^{\tau^*(c)\wedge\sigma}e^{-\rho t} g_c(c,D_t^{d_0})\ud t+ q_0^-e^{-\rho \sigma}\mathbf{1}_{\{\sigma<\tau^*(c)\}}
-q_0^+e^{-\rho \tau^*(c)}\mathbf{1}_{\{\tau^*(c)<\sigma\}}\right]\\
&\leq & \E\left[\int_0^{\tau\wedge\sigma}e^{-\rho t} g_c(c,D_t^{d_0})\ud t+ q_0^-\right]\\
&\leq & \E\left[\int_0^{\tau\wedge\sigma}e^{-\rho t} g_c(c,d)\ud t+ q_0^-\right]\\
&= &\frac{g_c(c,d)}{\rho}\,\E[1-e^{-\rho(\tau\wedge\sigma)}]+q^-_0.
\enqs
Note that $\sigma$ and $\tau$ are independent of $c$, and that $\tau\wedge\sigma>0$. So,  letting $c\rightarrow -\infty$ in the inequality above we get $\lim_{c\rightarrow-\infty}  {v_c}(c,d)=-\infty$, which contradicts Proposition \ref{prop:vc}\,(3).\ep
\begin{Remark}
{\rm We notice that items 1 and 2 of Corollary \ref{cor:fin} above hold, respectively, when $q_0^+<\infty$ and $q_0^-<\infty$, which is an assumption we are doing throughout the paper. However, also referring to Remark \ref{remrem}\,(2), we point out that in the case one consider, e.g., $q_0^-=\infty$ (irreversible investment), one has immediately $\hat{c}_-\equiv \infty$, so Corollary \ref{cor:fin} does not hold anymore. 
}
\end{Remark}
 \setlength{\unitlength}{2.5cm}
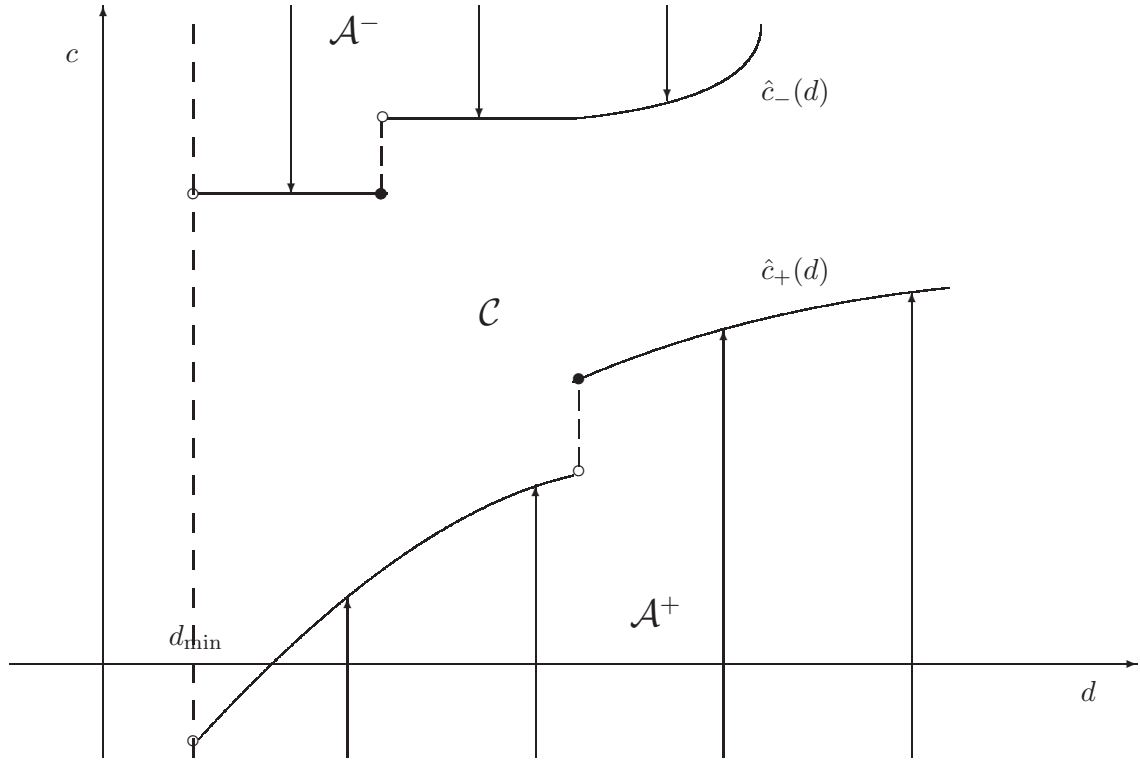
\begin{picture}(6,5)(-3,-2)
\put(-1.5,2.5){\vector(0,-1){1}}
\put(-0.5,2.5){\vector(0,-1){0.6}}
\put(.5,2.5){\vector(0,-1){0.5}}
\put(-3,-1){\vector(1,0){6}}
\put(-1.99,1.5){\line(1,0){1}}
\put(-1.99,1.501){\line(1,0){1}}
\put(-1.99,1.502){\line(1,0){1}}
\put(2.7,-1.2){$d$}
\put(-2.06,1.46){$\circ$}
\put(-1.06,1.46){$\bullet$}
\put(-.99,1.9){\line(1,0){1}}
\put(-.99,1.901){\line(1,0){1}}
\put(-.99,1.902){\line(1,0){1}}
\put(-.99,1.906){\line(1,0){1}}
\put(-1.05,1.87){$\circ$}
\put(-1.06,1.46){$\bullet$}
\put(-2.15,-0.9){$d_{\min}$}
\put(-.5,.8){$\mbox{\Large{$\mathcal{C}$}}$}
\put(.3,-.8){$\mbox{\Large{$\mathcal{A}^+$}}$}
\put(-1.3,2.3){$\mbox{\Large{$\mathcal{A}^-$}}$}
\put(-.01,-.012){$\circ$}
\put(-2.06,-1.45){$\circ$}
\put(-.01,.48){$\bullet$}
\put(-2.5,-1.5){\vector(0,1){4}}
\put(-1.2,-1.5){\vector(0,1){0.85}}
\put(-.2,-1.5){\vector(0,1){1.45}}
\put(.8,-1.5){\vector(0,1){2.28}}
\put(1.8,-1.5){\vector(0,1){2.48}}
\put(1,1.05){${\hat{c}_+}(d)$}
\put(1,2){${\hat{c}_-}(d)$}
\put(0.03,0.05){\line(0,1){.1}}
\put(0.03,0.20){\line(0,1){.1}}
\put(0.03,0.35){\line(0,1){.1}}
\put(-1.02,1.50){\line(0,1){.1}}
\put(-1.02,1.65){\line(0,1){.1}}
\put(-1.02,1.80){\line(0,1){.08}}
\put(-2.7,2.2)
{$c$}
\multiput(-2.02,-1.5)(0,.2){20}{\line(0,1){.1}}
\qbezier(0,.5)(0.8853,0.8853)
(2,1)
\qbezier(0,0)(-0.8853,-.1953)
(-1.99,-1.4)
\qbezier(0,.5)(0.8853,0.8853)
(2,1)
\qbezier(0,0)(-0.8853,-.1953)
(-1.99,-1.4)
\qbezier(0,.5)(0.8853,0.8853)
(2,1)
\qbezier(0,0)(-0.8853,-.1953)
(-1.99,-1.4)
\qbezier(0.0,1.9)(1,2)
(1,2.4)
\qbezier(0.0,1.9)(1,2)
(1,2.4)
\qbezier(0.0,1.9)(1,2)
(1,2.4)\qbezier(0.0,1.9)(1,2)
(1,2.4)
\qbezier(0.0,1.9)(1,2)
(1,2.4)
\qbezier(0.0,1.9)(1,2)
(1,2.4)
\end{picture}

\section{Quadratic cost:  smooth fit and boundaries' characterization}\label{sec:quadr}

\setcounter{equation}{0}
\setcounter{Assumption}{0} \setcounter{Theorem}{0}
\setcounter{Proposition}{0} \setcounter{Corollary}{0}
\setcounter{Lemma}{0} \setcounter{Definition}{0}
\setcounter{Remark}{0}

Theorem  \ref{prop:structure} and the continuity of $v_c$ in $\Sc$ yield  some optimality conditions. Indeed, we should have 
 \begin{equation}\label{smoothbase}
 \begin{cases}
A'(c)\psi(d)+B(c)\varphi (d)+ {{\hat{V}}}_c(c,d) \ \ = \ \ -{q_0^+},\ \ \ \ \ \ \forall \, (c,d)\in \partial\mathcal{C}^+,\\\\
A'(c)\psi(d)+B(c)\varphi (d)+ {{\hat{V}}}_c(c,d) \ \ = \ \ {q_0^-}, \  \ \ \ \ \ \ \ \forall \, (c,d)\in\partial \mathcal{C}^-.
\end{cases}
 \end{equation}
 
 It is clear that one cannot expect that the  conditions above provide a way either to find the value function or the optimal boundaries $\partial^{\pm}\mathcal{C}$ (e.g., in terms of the functions ${\hat{c}_{\pm}}$), as,  read at $(\hat{c}_\pm(d),d)$, they would relate four unknown functions $A,B, \hat{c}_{\pm}$ by two equations.
Other optimality conditions are needed and should be derived from some other suitable smoothness property of the value function at the optimal boundaries $\partial^\pm \Cc$.
To this end, we notice  by Theorem \ref{prop:structure} that
\begin{equation}\label{sf}
\frac{\partial}{\partial d}\, v_{c}(c,d)=0 \ \ \mbox{in} \ {\mathcal{A}^{\pm}}.
\end{equation}
Therefore, a requirement of a smooth fit condition of   
the second order mixed derivative of $v$ at the optimal boundaries would imply
\beq\label{sosf}
\lim_{(c,d)\rightarrow ({c}_0,{d}_0)} v_{cd}(c,d) \ \ = \ \ 0, \ \ \ \ \forall\, ({c}_0,{d}_0)\in \partial^\pm\mathcal{C}.
\enq
This is what we are going to prove in the next subsection under further assumptions on $g$.

\subsection{The smooth fit-principle}
The purpose  of the present subsection is indeed to prove  \eqref{sosf}. However,
we need to further specify our assumptions, restricting to  the quadratic cost case:  
 \beq \label{ggg}
g(c,d) &=& \frac{1}{2}(c^2- 2\beta_0(d)c+\alpha_0(d)),
\enq
where $\alpha_0, \beta_0$ are continuous functions. From now on, we assume that $g$ has the structure \eqref{ggg} and we do not repeat this assumption in the statements of the results.
We assume that  the functions $\alpha_0,\beta_0$ are continuous and that  $\beta_0$ is nondecreasing, so that Assumption \ref{ass:cost} holds true, and we denote 
\beq \label{alphabeta}
\alpha (d) \; := \; \E\Big[ \int_0^{\infty}e^{-\rho t} \alpha_0(D_t^d)  \ud t \Big],& &   \beta(d) \; := \; \E\Big[ \int_0^{\infty}e^{-\rho t} \beta_0(D_t^d) \ud t\Big],
\enq
noting that $\alpha,\beta\in C^2(\Oc;\mathbb{R})$ as the diffusion $D$ is nondegenerate.
The function ${\hat{V}}$ is  written in this case as: 
\beq \label{hatVqua}
{\hat{V}}(c,d) &=&  \frac{1}{2}\Big(\frac{1}{\rho}c^2-2\beta(d)c+\alpha(d)\Big).
\enq
Given a  function $\varphi\in C(\mathbb{R};\mathbb{R})$, let us denote  
\beqs
[\Delta^2\varphi](x;\varepsilon) &:=&  \frac{1}{\varepsilon^2} [\varphi(x+\varepsilon)+\varphi(x-\varepsilon)-2\varphi(x)], \ \ \ x\in\mathbb{R}, \ \varepsilon>0.
\enqs
The following Lemma, which relies  on  assumption \eqref{ggg},  enables us  to obtain  further regularity of the value function with respect to $c$ (Corollary \ref{cor:bene}), which is crucial to prove then \eqref{sosf}.

\begin{Lemma}\label{prop:secondorder}
We have for every $(c,d)\in\Sc,\  \varepsilon>0$,
\beqs
0 \ \  \leq \ \ [\Delta^2 v(\cdot,d)](c;\varepsilon) \ \ \leq  \ \  \frac{1}{\rho}\,.
\enqs
\end{Lemma}
\noindent \textbf{Proof.}
The estimate from below is a straightforward consequence of the convexity of $v$ 
with respect to $c$. Let us prove the estimate from above. Let $(c,d)\in{\mathcal{S}}$, $\varepsilon>0$, and $I\in\mathcal{I}$. 
By using the fact that  $g_{cc}\equiv 1$ under \reff{ggg}, we have
\beq
& & \frac{1}{\varepsilon^2}\left[G(c+\varepsilon,d;I)+G(c-\varepsilon,d;I)-2 G(c,d;I)\right] \label{qdr} \\
&=& \mathbb{E}\Big[\int_0^{\infty}e^{-\rho t} \Big[\frac{1}{\varepsilon^2}\big(g(C^{c+\varepsilon,I}_t,D_t^d)+g(C^{c-\varepsilon,I}_t,D_t^d)-2g(C^{c,I}_t,D_t^d)\big)\Big) \ud t\Big] \; = \;  \frac{1}{\rho}. \nonumber
\enq
%
Since
\beqs
v(c+\varepsilon,d)+v(c-\varepsilon,d)-2 G(c,d;I) &\leq& G(c+\varepsilon,d;I)+G(c-\varepsilon,d;I)-2 G(c,d;I),
\enqs
we get from  \eqref{qdr}: 
\beqs
\frac{1}{\varepsilon^2}\left[v(c+\varepsilon,d)+v(c-\varepsilon,d)-2 G(c,d;I)\right] &\leq&  \frac{1}{\rho}, \ \ \ \ \ \forall I\in\Ic.
\enqs
Taking the supremum over $I\in\mathcal{I}$, this proves  the required upper-estimate.
\ep

\vspace{2mm}

Lemma \ref{prop:secondorder} implies that  $v_c(\cdot,d)$ is  Lipschitz continuous for each $d\in\mathcal{O}$. Together with \eqref{AB}-\eqref{AB2} and \eqref{hatVqua}, we  immediately  get the following regularity result.

\begin{Corollary}\label{cor:bene}
The derivative functions $A',B':(\underline{c}_+,\bar{c}_-)\rightarrow \mathbb{R}$, where $A,B$ are the functions defined in Theorem \ref{prop:structure}, are locally Lipschitz. 
In other terms $A,B\in W_{loc}^{2,\infty}((\underline{c}_+,\bar{c}_-);\mathbb{R})$.

(This property holds eventually up to $\underline{c}_+,\bar{c}_-$, when $A,B$ can be  extended, respectively, to $C^1$ functions up to $\underline{c}_+,\bar{c}_-$, according to the conditions of Theorem  \ref{prop:structure} which allow these extensions.)
\end{Corollary}

\vspace{2mm}

We are now able to prove the second order smooth-fit  result on the value function.

\begin{Proposition}\label{prop:sosf}
The relation \eqref{sosf} hold true.
\end{Proposition}
\noindent \textbf{Proof.} Since $v_{cd} = 0$ in $\mathcal{A}^\pm$,  the claim is equivalent to prove that 
\beq\label{sosf2}
\lim_{\stackrel{(c,d)\rightarrow ({c}_0,{d}_0)}{(c,d)\in\Cc}} v_{cd}(c,d) \ \ = \ \ 0, \ \ \ \ \forall\, ({c}_0,{d}_0)\in \partial^\pm\mathcal{C}.
\enq
We shall prove \eqref{sosf2} for the lower boundary $\partial^+ {\Cc}$; the claim concerning the upper boundary $\partial^-{\Cc}$ can be proved in the same way. Letting $(c_0,d_0)\in \partial^+\Cc$ we distinguish three cases.

1. Suppose that $c_0=\hat{c}_+(d)>\underline c_+$. Let us  consider the function  on $\mathcal{D}:=(\underline{c}_+,\bar{c}_-)\times \mathcal{O}$ 
\beqs
\bar{v}(c,d) &:=&   A(c)\psi(d)+B(c)\varphi(d)+{\hat{V}}(c,d), \ \ \ \ (c,d)\in\Dc.
\enqs
  By Theorem  \ref{prop:structure} and \eqref{hatVqua}, we have that $\bar{v}\in C^{1,2}(\mathcal{D};\R)$, and that $\bar{v}_{cd}$ exists and is continuous in $\mathcal{D}$. Since $\bar{v}=v$ in $ \bar{\mathcal{C}}\cap \mathcal{D},$ by monotonicity of $v_c(c,\cdot)$, we have 
  \beq\label{iij}
  \bar{v}_{cd}&\leq& 0 \ \ \ \mbox{in}\  \Cc.
  \enq
Clearly   \eqref{sosf2} is equivalent to 
\beq\label{pkjh}
\lim_{\stackrel{(c,d)\rightarrow ({c}_0,{d}_0)}{(c,d)\in\Cc}} \bar{v}_{cd}(c,d) \ \ = \ \ 0, \ \ \ \ \forall\, ({c}_0,{d}_0)\in \partial^+\mathcal{C}.
\enq
By continuity of $\bar{v}_{cd}$, the limit above exists and coincides with $\bar{v}_{cd}(c_0,d_0)$. Taking into account \eqref{iij}, suppose  by contradiction  that
\beq \label{psd}
\bar{v}_{cd}(c_0,d_0) &<& 0.
\enq
Then, by continuity of $\bar{v}_{cd}$, we may find 
$\varepsilon>0$, $\delta_1>0$, $\delta_2>0$  such that   
\beq \label{psd2}
\bar{v}_{cd}(c,d) &\leq&  -\varepsilon, \ \ \ \forall (c,d)\in (c_0-\delta_1,c_0-\delta_1)\times (d_0-\delta_2,d_0+\delta_2)\subset \mathcal{D}.
\enq
Since  $\bar{v}_{c}(c_0,d_0)=-q_0^+$, 
due to \eqref{psd} and to Corollary \ref{cor:bene}, we can apply Implicit Function Theorem  in a generalized form,
 stating that  there exists  $\hat{d}_+'$  in Sobolev sense in the interval $(c_0-\delta_1,c_0+\delta_1)$, and, assuming without loss of generality that 
 \beq\label{wlog}
 d_0-\delta_2 \ = \ \hat{d}_+(c_0-\delta_1), \ \ \ \  d_0+\delta_2 \  = \ \hat{d}_+(c_0+\delta_1),
 \enq
that  it holds, by Corollary \ref{cor:bene} 
and \eqref{psd2}
\beq\label{trw}\nonumber
\hat{d}_+' (\cdot) \ \  =  \  \  -\frac{\bar{v}_{cc}(\cdot,\hat{d}_+(\cdot))}{\bar{v}_{cd}(\cdot,\hat{d}_+(\cdot))} & \leq & 
M_\varepsilon\  \ <\ \ \infty,\  \ \  \mbox{a.e. in } \ (c_0-\delta_1,c_0+\delta_1).
\enq 
Let us now assume, without loss of generality (recall that $\hat{c}_+$ is right-continuous), that $\hat{c}_+$ is continuous on $[d_0,d_0+\delta_2)$. 
Then, combining with \eqref{wlog}-\eqref{trw}, we see that $\hat{c}_+$ is strictly increasing on $[d_0,d_0+\delta_2)$, there exists the inverse $\hat{c}_+^{-1}$ on $[c_0,c_0+\delta_1)$, it coincides with $\hat{d}_+$, and $\hat{d}_+$ is continuous and strictly increasing as  well on $[c_0,c_0+\delta_1)$. It follows that
\beq\label{ddf}
\hat{d}_+'&>& 0, \ \ \ \mbox{a.e. in} \ \ [c_ 0,c_0+\delta_1). 
\enq
Let   $\Yc$ be the set of differentiability points of $\hat{d}_+$ in $[d_0,d_0+\delta)$ where $0<\hat{d}_+'<M_\varepsilon$. Then, taking into account \eqref{trw}-\eqref{ddf}, we see that $\mathcal{Y}$ has full measure in $[c_0,c_0+\delta_1)$. Consequently $\hat{d}(\mathcal{Y})$ is dense in $[d_0,d_0+\delta_2)$, $\hat{c}_+'$ exists in $\hat{d}_+(\Yc)$,   and
\beq\label{kkj}
\hat{c}'_+ \ \ \in \ \ [1/M_\varepsilon,\infty), \ \ \ \ \ \mbox{in}\ \ \hat{d}_+(\Yc).
\enq  
Let us now consider the function $d \in [d_0,d_0+\delta_2)$ $\mapsto$ $v(c_0,d).$ 
Since ${\hat{c}_+}$ is nondecreasing  in $[d_0,d_0+\delta_2)$ (actually we have shown strictly increasing), the segment $\{(c_0,d) \ | \ d\in [d_0,d_0+\delta_2)\}$ is contained in $\mathcal{A}^+$. Hence,  Theorem \ref{prop:structure} yields
\begin{equation}\label{VVV}
v(c_0,d)\ \ =\ \ -{q_0^+}c_0+z_+(d), \ \ \ \ \ \forall d\in [d_0,d_0+\delta).
\end{equation}
Applying the chain rule at the points of $\hat{d}_+(\mathcal{Y})$ to 
$$
[d_0,d_0+\delta)\ \rightarrow\ \mathbb{R},  \;\;\; \ \ \  d\ \mapsto \ z_+(d) \; = \; v({\hat{c}_+}(d),d)+{q_0^+}{\hat{c}_+}(d) \;= \; \bar{v}({\hat{c}_+}(d),d)+{q_0^+}{\hat{c}_+}(d),$$
 we see that the function $z_+$ is differentiable at the points of  $\hat{d}_+(\mathcal{Y})$ and 
\beqs  
z_+'(d) &=& \bar{v}_{c}({\hat{c}_+}(d),d)\hat{c}_+'(d)+\bar{v}_d({\hat{c}_+}(d),d)+{q_0^+}\hat{c}_+'(d), \ \ \ \forall  d\in\hat{d}_+(\Yc).
\enqs
By definition of ${\hat{c}_+}$, we have $\bar{v}_{c}({\hat{c}_+}(d),d)=v_c({\hat{c}_+}(d),d)=-{q_0^+}$ for every $d\in\mathcal{O}$, and so
\beqs \label{VVV1}
z'_+(d) \; = \; v_d({\hat{c}_+}(d),d),\ \ \ \forall d\in\hat{d}_+(\mathcal{Y}).
\enqs
Together with  \eqref{VVV}, this shows   the existence of $v_d(c_0,d)$ for each $d\in\hat{d}_+(\mathcal{Y})$ and the equality
\begin{equation} \label{snd}
v_d(c_0,d) \; = \; z_+'(d) \; = \; \bar{v}_d({\hat{c}_+}(d),d),\ \ \ \forall d\in\mathcal{Y}.
\end{equation}
On the other hand,  by using again the chain rule, we can  get from \eqref{snd} the existence of $v_{dd}(c_0,d)$ for each $d\in\hat{d}_+({\Yc})$ and the equality
\beq\label{z''}
v_{dd}(c_0,d) \; = \; z_+''(d) \; = \; \bar{v}_{dd}({\hat{c}_+}(d),d)+\bar{v}_{cd}({\hat{c}_+}(d),d)\,\hat{c}_+'(d),\ \ \ \forall  d\in\hat{d}_+(\mathcal{Y}).
\enq
Therefore, from  \eqref{psd2}, \eqref{kkj},  and \eqref{z''}, we get
\begin{equation}\label{VDD}
v_{dd}(c_0,d) \ \ \leq \ \    \bar{v}_{dd}({\hat{c}_+}(d),d)-\varepsilon/M_\varepsilon,\ \ \ \forall  d\in\hat{d}_+(\mathcal{Y}).
\end{equation}
Now the viscosity subsolution property of $v$, and \eqref{VVV}, \eqref{snd}, \eqref{VDD} yield
\beq\label{mza}
g(c_0,d) &\geq & \rho v(c_0,d)-\mu(d)v_{d}(c_0,d)-\frac{1}{2}\sigma(d)^2v_{dd}(c_0,d) \\
&=&\nonumber
\rho v(c_0,d)-\mu(d)\bar{v}_{d}({\hat{c}_+}(d),d)-\frac{1}{2}\sigma(d)^2[\bar{v}_{dd}({\hat{c}_+}(d),d)-\varepsilon/M_\varepsilon], \ \ \  \forall  d\in\hat{d}_+(\mathcal{Y}).
\enq
Taking a sequence $(\alpha_n)\subset\hat{d}_+(\mathcal{Y})$ such that $\alpha_n\downarrow d_0$ (this can be done since $\hat{d}_+(\Yc)$ is dense in $[d_0,d_0+\delta_2)$) and passing to the limit in  \eqref{mza} evaluated at $d=\alpha _n$ we obtain by continuity of $\hat{c}_+$ in $[d_0,d_0+\delta_2)$, continuity of $g$ in $\Sc$, and since $\bar{v}\in C^{1,2}(\Dc,\R)$ anfd $\bar{v}=v$ in $\bar{\Cc}$, 
\begin{equation}\label{mza1}
\rho \bar{v}(c_0,d_0)-\mu(d_0)\bar{v}_{d}(c_0,d_0)-\frac{1}{2}\sigma(d_0)^2[\bar{v}_{dd}(c_0,d_0)-\varepsilon/M_\varepsilon] \ \ \leq \ \   g(c_0,d_0).
\end{equation}
On the other hand, recall that  $\mathcal{L}\bar{v}=\mathcal{L}v=g$ on $\mathcal{C}$. 
Therefore, since $v\in C^{1,2}(\mathcal{D};\mathbb{R})$ and since $(c_0,d_0)\in\bar\Cc$, by continuity  we must also have 
\beqs
\rho \bar{v}(c_0,d_0)-\mu(d_0)\bar{v}_{d}(c_0,d_0)-\frac{1}{2}\sigma(d_0)^2\bar{v}_{dd}(c_0,d_0) &=& g(c_0,d_0),
\enqs
which is in contradiction with  \eqref{mza1}  as  $\sigma^2(d_0)>0$, and the claim is proved in this case.

2. Consider now the case $c_0=\hat{c}_+(d_0)=\underline{c}_+$. In this case we can construct the function $\hat{v}$ in $\mathcal{D}:=(\underline{c}_+-\varepsilon,\bar{c}_-)\times \mathcal{O}$ for some $\varepsilon>0$ by  using the extension part of Corollary \ref{cor:bene}, and repeat the argument of the previous case.  

3. Consider now the last possible case, i.e. $d_0=\hat{d}_+(c_0)$ and $c_0<\hat{c}_+(d_0)$, noting that $\hat{c}_+(d_0)<\infty$ (see Proposition \ref{frontiera}\,(1)). In this case the segment $K:=\{ (c,d_0) \ | \ c\in [c_0,\hat{c}_+(d_0)]\}$ is contained in $\partial ^+\Cc$.
Define the function $\bar{v}$ as in item 1. We then have $\bar{v}_c = v_c = -q_0^+$ in $K$. Hence 
\beq\label{bbnx}
-q_0^+-\bar{v}_c(c,d)  & = & \bar{v}_c(c,d_0)-\bar{v}_c(c,d)\nonumber \\
&=& \int_{d}^{d_0} \bar{v}_{cd}(c,\xi)\ud \xi, \ \ \ \ \  \forall c\in [c_0,\hat{c}_+(d_0)], \ \forall d\leq d_0, 
\enq
Taking into account Corollary \ref{cor:bene} and differentiating \eqref{bbnx} with respect to $c$ we get (the derivatives $A'',B''$ must be intended  in Sobolev sense)
\beq\label{bbnx1}
-\bar{v}_{cc}(c,d)  
&=& \int_{d}^{d_0} \bar{v}_{cdc}(c,\xi)\ud \xi, \ \ \ \ \ \   \mbox{a.e.} \ \ (c,d)\in [c_0,\hat{c}_+(d_0)]\times(\hat{d}_-(c),d_0].
\enq
Since $v_{cc}\geq 0$, hence $\bar{v}_{cc}\geq 0$ (in Sobolev sense),  from \eqref{bbnx1} we get 
\beq\label{bbnx2}
0&\geq & \int_{d}^{d_0} \bar{v}_{cdc}(c,\xi)\ud \xi, \ \ \ \ \ \   \mbox{a.e.} \ \ (c,d)\in [c_0,\hat{c}_+(d_0)]\times(\hat{d}_-(c),d_0],
\enq
from which, taking into account \eqref{hatVqua}, we deduce that actually 
\beqs
A''(c)\psi'(d)+B''(c)\varphi'(d)&\leq &0, \ \ \ \ \ \   \mbox{a.e. in} \ \ [c_0,\hat{c}_+(d_0)]\times(\hat{d}_-(c),d_0],
\enqs
Then, since $\psi',\varphi'$ are continuous, we deduce that 
\beqs
A''(c)\psi'(d_0)+B''(c)\varphi'(d_0)&\leq &0, \ \ \ \ \ \   \mbox{a.e. in} \ \ [c_0,\hat{c}_+(d_0)].
\enqs
Hence, $\bar{v}_{cd}(\cdot,d_0)$ is nonincreasing with respect to $c$ in $[c_0,\hat{c}_+(d_0)]$. Then,
assuming now, as in item 1, by contradiction \eqref{psd}, we also must have $\bar{v}_{cd}(\hat{c}_+(d_0),d_0)<0$. So we are now reduced to the contradiction assumption of item 1, we can apply the argument of that item and get the contradiction, so the claim.   
\ep

\begin{Remark} 
{\rm  In \cite{MZ}, a similar  smooth-fit principle \reff{sosf} is derived a posteriori in the particular case where the 
state process is a geometric Brownian motion, so that an explicit smooth solution can be obtained, and then shown to be the equal to the value function by a verification approach.  In the general  diffusion case for demand and when the cost function is quadratic, we prove directly the smooth-fit principle \reff{sosf} by a viscosity solutions approach.  
}
\end{Remark}

\subsection{Characterization of the optimal boundaries}
Proposition \ref{prop:sosf}  can be used to  add other necessary optimality conditions to \eqref{smoothbase}:  indeed, by  \eqref{structureclosure}, the relation  \eqref{sosf} yields
\beq\label{opt2}
  A'(c)\psi'(d)+B'(c)\varphi'(d)+\hat{V}_{cd}(c,d) \ \ =\ \  0, \ \ \ \forall\, (c,d)\in \partial \Cc
 \enq
 We want to use  the optimality conditions  \eqref{smoothbase} and \eqref{opt2}  to characterize  the optimal boundaries $\partial{\mathcal{C}}^{\pm}$. First, we rewrite such conditions. (The proofs of the next two propositions follow the line of \cite{A2} and also, in some parts,  of \cite{MZ}.) 
   \begin{Proposition}\label{prop:ex}
Let $c\in \R$ and let $d_+,d_-\in\Oc$ be such that $(c,d_-)\in\partial^-\Cc$, $(c,d_+)\in\partial^+\Cc$. Then
    \beq\label{optcond}
\begin{cases}
{\displaystyle{\int_{{d}_-}^{{d}_+}\psi(\xi)g_c(c,\xi)m'(\xi)\ud \xi\ +\ q_0^-\frac{\psi'({d}_-)}{S'({d}_-)}\ +\ q_0^+\frac{\psi'({d}_+)}{S'({d}_+)}\ \ =\ \ 0,}}\\\\
{\displaystyle{\int_{{d}_-}^{{d}_+}\varphi(\xi)g_c(c,\xi)m'(\xi)\ud \xi\ + \ q_0^-\frac{\varphi'({d}_-)}{S'({d}_-)}\ +\  q_0^+\frac{\varphi'({d_+})}{S'({d}_+)}\ \ =\ \ 0.}}
\end{cases}
  \enq
   \end{Proposition}
  
 \noindent \textbf{Proof.}
Let $c, d_\pm$ be as in the statement. 
The conditions \eqref{smoothbase} computed respectively at $(c,d_+)$ and $(c,d_-)$ yield
\beqs
\begin{cases}
A'(c)\psi(d_+)+B'(c)\varphi(d_+)+\hat{V}_{c}(c,d_+)\ \ =\ \ -q_0^+,\\
A'(c)\psi(d_-)+B'(c)\varphi(d_-)+\hat{V}_{c}(c,d_-)\ \ =\ \ q_0^-,\\
\end{cases}
\enqs
from which we get
\beq\label{eds}
\begin{cases}
A'(c)\ \ =\  \ \frac{\varphi(d_-)(-\hat{V}_{c}(c,d_+)-q_0^+)-\varphi(d_+)(q_0^--\hat{V}_{c}(c,d_-))}{\psi(d_+)\varphi(d_-)-\varphi(d_+)\psi(d_-)},\\\\
B'(c)\ \ =\  \ \frac{\psi(d_+)(q_0^--\hat{V}_c(c,d_-))-\psi(d_-)(-q_0^+-\hat{V}_c(c,d_+))}{\psi(d_+)\varphi(d_-)-\varphi(d_+)\psi(d_-)}.
\end{cases}
\enq
By Theorem \ref{prop:structure}
\beq\label{hgf}
v_c(c,d)&=&A'(c)\psi(d)+B'(c)\varphi(d)+\hat{V}_{c}(c,d),\ \ \forall d\in [d_-,d_+].
\enq
So, plugging \eqref{eds}  into \eqref{hgf},  
we get
\beq
v_c(c,d)\;=\; \frac{\tilde{\varphi}(d)}{\tilde{\varphi}(d_-)}(q_0^--\hat{V}_{c}(c,d_-))+\frac{\tilde{\psi}(d)}{\tilde{\psi}(d_+)}(-q_0^+-\hat{V}_{c}(c,d_+))+\hat{V}_{c}(c,d), \   \forall d\in [d_-,d_+],
 \enq
where
\beq
\tilde{\varphi}(d)&:=& \varphi(d)-\frac{\varphi(d_+)}{\psi(d_+)}\psi(d),\ \ \ \ \tilde{\psi}(d)\ \ :=\ \  \psi(d)-\frac{\psi(d_-)}{\varphi(d_-)}\varphi(d).
\enq
 Hence
 \beq\label{vcc}
 v_{cd}(c,d) \;=\;\frac{\tilde{\varphi}'(d)}{\tilde{\varphi}(d_-)}(q_0^--\hat{V}_{c}(c,d_-))+\frac{\tilde{\psi}'(d)}{\tilde{\psi}(d_+)}(-q_0^+-\hat{V}_{c}(c,d_+))+ \hat{V}_{cd}(c,d), \, \forall d\in [d_-,d_+].
  \enq
Now \eqref{opt2} yields $v_{cd}(c,d_-)=v_{cd}(c,d_+)=0$. Imposing these conditions into \eqref{vcc},  we get
\beq\label{jjh}
\begin{cases}
q_0^--\hat{V}_{c}(c,d_-)\ \ =\ \ \frac{-\hat{V}_{cd}(c,d_-) \tilde{\psi}'(d_+)\tilde{\varphi}(d_-)+ \hat{V}_{cd}(c,d_+)\tilde{\psi}'(d_-)\tilde{\varphi}(d_-)}{\tilde{\varphi}'(d_-)\tilde{\psi}'(d_+)-\tilde{\psi}'(d_-)\tilde{\varphi}'(d_+)},\\
-q_0^+-\hat{V}_{c}(c,d_+)\ \ =\ \ \frac{-\hat{V}_{cd}(c,d_+) \tilde{\varphi}'(d_-)\tilde{\psi}(d_+)+ \hat{V}_{cd}(c,d_-)\tilde{\varphi}'(d_+)\tilde{\psi}(d_+)}{\tilde{\varphi}'(d_-)\tilde{\psi}'(d_+)-\tilde{\psi}'(d_-)\tilde{\varphi}'(d_+)}.
\end{cases}
\enq
Simple computations yield
\beqs
\tilde{\varphi}'(d_-)\tilde{\psi}'(d_+)-\tilde{\psi}'(d_-)\tilde{\varphi}'(d_+)&=& (\varphi'(d_-)\psi'(d_+)-\varphi'(d_+)\psi'(d_-))(\varphi(d_-)\psi(d_+)-\varphi(d_+)\psi(d_-)),\\
\tilde{\psi}'(d_+)\tilde{\varphi}(d_-)&=& \frac{(\psi'(d_+)\varphi(d_-)-\psi(d_-)\varphi'(d_+))(\varphi(d_-)\psi(d_+)-\varphi(d_+)\psi(d_-))}{\psi(d_+)\varphi(d_-)},\\
\tilde{\psi}'(d_-)\tilde{\varphi}(d_-)&=& \frac{(\psi'(d_-)\varphi(d_-)-\psi(d_-)\varphi'(d_-))(\varphi(d_-)\psi(d_+)-\varphi(d_+)\psi(d_-))}{\psi(d_+)\varphi(d_-)},\\
\tilde{\varphi}'(d_-)\tilde{\psi}(d_+)&=& \frac{(\varphi'(d_-)\psi(d_+)-\varphi(d_+)\psi'(d_-))(\varphi(d_-)\psi(d_+)-\varphi(d_+)\psi(d_-))}{\psi(d_+)\varphi(d_-)},\\
\tilde{\varphi}'(d_+)\tilde{\psi}(d_+)&=& \frac{(\varphi'(d_+)\psi(d_+)-\varphi(d_+)\psi'(d_+))(\varphi(d_-)\psi(d_+)-\varphi(d_+)\psi(d_-))}{\psi(d_+)\varphi(d_-)}.
\enqs
Plugging these expressions into \eqref{jjh} we get
\beq\label{paswe}
\begin{cases}
q_0^--\hat{V}_{c}(c,d_-)\ \ =\ \ \frac{-\hat{V}_{cd}(c,d_-) (\psi'(d_+)\varphi(d_-)-\psi(d_-)\varphi'(d_+))+ \hat{V}_{cd}(c,d_+)(\psi'(d_-)\varphi(d_-)-\psi(d_-)\varphi'(d_-)}{{\varphi}'(d_-){\psi}'(d_+)-{\psi}'(d_-){\varphi}'(d_+)},\\
-q_0^+-\hat{V}_{c}(c,d_+)\ \ =\ \ \frac{-\hat{V}_{cd}(c,d_+)(\varphi'(d_-)\psi(d_+)-\varphi(d_+)\psi'(d_-))+ \hat{V}_{cd}(c,-)(\varphi'(d_+)\psi(d_+)-\varphi(d_+)\psi'(d_+))}{{\varphi}'(d_-){\psi}'(d_+)-{\psi}'(d_-){\varphi}'(d_+)}.
\end{cases}
\enq
Using the representations \eqref{Vc}-\eqref{Vcd} in \eqref{paswe}, we get after long computations
\begin{eqnarray*}
-q_0^+({\varphi}'(d_-){\psi}'(d_+)-{\psi}'(d_-){\varphi}'(d_+)) & =&  \varphi'(d_-)S'(d_+)\int_{d_-}^{d_+} \psi(\xi)g_c(c,\xi)m'(\xi)\ud\xi\\
&& -\psi'(d_-)S'(d_+)\int_{d_-}^{d_+}\varphi(\xi)g_c(c,\xi)m'(\xi)\ud \xi,
\end{eqnarray*}
\begin{eqnarray*}
q_0^-({\varphi}'(d_-){\psi}'(d_+)-{\psi}'(d_-){\varphi}'(d_+))& =&  \varphi'(d_+)S'(d_-)\int_{d_-}^{d_+} \psi(\xi)g_c(c,\xi)m'(\xi)\ud\xi
\\
&&-\psi'(d_+)S'(d_-)\int_{d_-}^{d_+}\varphi(\xi)g_c(c,\xi)m'(\xi)\ud \xi,
\end{eqnarray*}
from which we finally see that the couple $(d_-,d_+)\in\Oc\times\Oc$ satisfies  \eqref{optcond}.
  \ep\\\\
Let us denote 
$$\underline{c}_{+,g} \ :=\ \inf_{\mathcal{O}} \hat{c}_{+,g}, \ \ \ \ \underline{c}_{-,g} \ :=\ \inf_{\mathcal{O}} \hat{c}_{-,g}, \ \ \  \ \bar{c}_{+,g} \ :=\ \sup_{\mathcal{O}} \hat{c}_{+,g},\ \  \ \   \bar{c}_{-,g} \ :=\ \sup_{\mathcal{O}} \hat{c}_{-,g}.$$
For all $c\in\R$ denote 
$$
d^*_+(c)\ \ := \ \ \inf\,\{\xi\in\Oc \ |\  g_c(c,\xi)<-\rho q_0^+\}, \ \ \ d^*_-(c)\ \ :=\ \ \sup\,\{\xi\in\Oc \ |\  g_c(c,\xi)>\rho q_0^-\}. 
$$ 
with the convention $\sup \emptyset=d_{\min}$, $\inf \emptyset =d_{\max}$. Then clearly we have $d^*_+(c)<d^*_-(c)$ for every $c\in \R$, and 
$d^*_+(c),d^*_-(c)\in \Oc$ if and only if $c \in (\underline{c}_{-,g}, \bar{c}_{+,g})$.
\begin{Proposition}\label{prop:exun}
Let $c\in \R $ and let $-\beta_0$ be strictly decreasing (so that $g_c(c,\cdot)=-\beta_0(\cdot)$ is strictly decreasing for every $c\in\R$). The couple of equations
  \beq\label{optcond0}
\begin{cases}
{\displaystyle{\int_{x}^{y}\psi(\xi)g_c(c,\xi)m'(\xi)\ud \xi\ + \ q_0^-\frac{\psi'(x)}{S'(x)}\ +\ q_0^+\frac{\psi'(y)}{S'(y)}}} \ \ =\ \ 0\\\\
{\displaystyle{\int_{x}^{y}\varphi(\xi)g_c(c,\xi)m'(\xi)\ud \xi\ + \ q_0^-\frac{\varphi'(x)}{S'(x)}\ +\ q_0^+\frac{\varphi'(y)}{S'(y)}\ \ =\ \ 0.}}
\end{cases}
\enq
admits a solution $(x^*(c),y^*(c))$ with $y^*(c)>x^*(c)$ if and only if  $c\in (\underline{c}_{-,g}, \bar{c}_{+,g})$ (note that the case $\underline{c}_{-,g}> \bar{c}_{+,g}$ may occur, and in this case this interval is considered  as empty). If this is the case, i.e.  $c\in (\underline{c}_{-,g}, \bar{c}_{+,g})$, then the  solution is unique and  belongs to 
$(d_{\min},d^*_-(c))\times (d^*_+(c),d_{\max})$.

Moreover $x^*, y^*$ are continuously differentiable in the interval $(\underline{c}_{-,g}, \bar{c}_{+,g})$ and have strictly positive derivatives.
\end{Proposition}
\noindent \textbf{Proof.}  
Fix $c\in\R$ and
  consider the functions in the couple of variables $(x,y)\in \Oc\times\Oc$ 
\beq\label{L1}
L_1(x,y;c) & :=& \int_{x}^y \psi(\xi)g_c(c,\xi)m'(\xi)\ud \xi\ +\ q_0^+ \frac{\psi'(y)}{S'(y)}\ +\ q_0^-\frac{\psi'(x)}{S'(x)},
\enq
\beq\label{L2}
L_2(x,y;c) & :=& \int_{x}^y \varphi(\xi)g_c(c,\xi)m'(\xi)\ud \xi\ +\ q_0^+\frac{\varphi'(y)}{S'(y)}\ +\ q_0^-\frac{\varphi'(x)}{S'(x)}.
\enq
The solvability of our system of equations corresponds then to the solvability of $L_1(x,y;c)=0$, $L_2(x,y;c)=0$ in $\Oc\times \Oc$ with $x<y$.
Using the representations (see, e.g., \cite[Ch.\,II]{BS})
\beq\label{repr}
\frac{\psi'(\cdot)}{S'(\cdot)}\ \ =\ \ \rho \int_{d_{\min}}^{\cdot}\psi(\xi)m'(\xi)\ud \xi,
\ \ \ \ \ \ \frac{\varphi'(\cdot)}{S'(\cdot)}\ \ =\ \ - \rho \int_{\cdot}^{d_{\max}}\varphi(\xi)m'(\xi)\ud \xi,
\enq
$L_1,L_2$ can be rewritten as 
\beqs
L_1(x,y;c)  &=& \int_{x}^y \psi(\xi)(g_c(c,\xi)+\rho q_0^+)m'(\xi)\ud \xi\ +\ (q_0^++q_0^-)\frac{\psi'(x)}{S'(x)},
\enqs 
\beqs
L_2(x,y;c)  &=& \int_{x}^y \varphi(\xi)(g_c(c,\xi)-\rho q_0^-)m'(\xi)\ud \xi\ +\ (q_0^++q_0^-)\frac{\varphi'(y)}{S'(y)},
\enqs 
or equivalently as
\beqs
L_1(x,y;c) & =& \int_{x}^y \psi(\xi)g_c(c,\xi)m'(\xi)\ud \xi\ +\ \rho q_0^+ \int_{d_{\min}}^y \psi(\xi)m'(\xi)\ud \xi\ + \ \rho q_0^- \int_{d_{\min}}^x \psi(\xi)m'(\xi)\ud \xi,
\enqs
\beqs
L_2(x,y;c) & =& \int_{x}^y \varphi(\xi)g_c(c,\xi)m'(\xi)\ud \xi\ \ -\ \rho q_0^+\int_y^{d_{\max}} \varphi(\xi)m'(\xi)\ud \xi \ -\ \rho q_0^-\int_x^{d_{\max}} \varphi(\xi)m'(\xi)\ud \xi,
\enqs
and the partial derivatives of $L_1,L_2$ with respect to $x,y$
are
$$\frac{\partial L_1}{\partial x} (x,y;c)\ \ =\ \ -\psi(x)(g_c(c,x)-\rho q_0^-)m'(x), \ \ \ \ \frac{\partial L_1}{\partial y} (x,y;c)\ \ =\ \ \psi(y)(g_c(c,y)+\rho q_0^+)m'(y),$$
$$
\frac{\partial L_2}{\partial x} (x,y;c)\ \ =\ \ -\varphi(x)(g_c(c,x)-\rho q_0^-)m'(x), \ \ \ \ \frac{\partial L_2}{\partial y} (x,y;c)\ \ =\ \ \varphi(y)(g_c(c,y)+\rho q_0^+)m'(y).
$$

Let us study the solvability of $L_1(x,\cdot;c)=0$ for given $x\in\Oc$. 
First of all we notice that $L_1(x,x;c)>0$ as $\psi'>0$, $S'>0$.
Taking into account that $g_c(c,\cdot)$ is strictly decreasing and continuous, we see that  the sign of $\frac{\partial L_1}{\partial y}(x,\cdot;c)$ is strictly positive in $(x,d_+^*(c))$ and strictly negative in $(d^*_+(c),d_{\max})$. Combined with the fact that 
$L_1(x,x;c)>0,$ this shows that there is at most one point $y^*(x;c)\in(x,d_{\max})$  solution to  $L_1(x,\cdot;c)=0$ and that $y^*(x;c)$ (if exists) must belong to $(d^*_+(c),d_{\max})$. 
Now we distinguish two cases.
\begin{itemize}
\item[-] If $c\geq \bar{c}_{+,g}$, then $g_c(c,\cdot)+\rho q_0^+\geq 0$ in $\Oc$. So the solution does not exist in this case.
\item[-] If  $c< \bar{c}_{+,g}$, take $\hat{y}(c)>d^*(c)$ such that $L_1(x,\hat{y}(c);c)>0$ (such $\hat{y}(c)$ exists by continuity), and observe that since $g_c(c,\cdot)$ is (strictly) decreasing, using \eqref{repr}, one has for every $y\geq \hat{y}(c)$
$$
\int_{\hat{y}(c)}^y\psi(\xi)m'(\xi)(g_c(c,\xi)+\rho q_0^+)\ud \xi\ \ \leq\ \ \frac{g_c(c,\hat{y})+\rho q_0^+}{\rho}\left(\frac{\psi'(y)}{S'(y)}-\frac{\psi'(\hat{y})}{S'(\hat{y})}\right),
$$
 therefore
 \beq\label{rrt}
 L_1(x,y;c)\ \ \leq \ \ L_1(x,\hat{y}(c);c)\ +\ \frac{g_c(c,\hat{y}(c))+\rho q_0^+}{\rho}\left(\frac{\psi'(y)}{S'(y)}-\frac{\psi'(\hat{y})}{S'(\hat{y})}\right).
 \enq
Now we notice that there exists $M_c>0$ such that $L_1(x,\hat{y}(c);c)\leq M_c$ for every $x\leq \hat{y}(c)$.  Indeed, 
$ \int_{d_{\min}}^{\hat{y}(c)} \psi(\xi)g_c(c,\xi)m'(\xi)\ud \xi$ is finite because of the finiteness of $\hat{V}_c$ and taking into account \eqref{Vcd};
$\int_{d_{\min}}^{\hat{y}(c)} \psi(\xi)m'(\xi)\ud \xi$ is finite because of \eqref{repr};  ${\psi'(x)}/{S'(x)}$ is bounded in $(d_{\min},\hat{y}(c)]$ because of \eqref{psiphi'}.
Now, 
since  $g_c(c,\hat{y}(c))+\rho q_0^+<0$ and since by \eqref{psiphi'} we have $\psi'(y)/S'(y)\rightarrow \infty$ as $y\rightarrow d_{\max}$, we see that the solution  $y^*(x;c)$ to $L_1(x,\cdot;c)=0$ exists in the interval $(\hat{y}(c),d_{\max}-\varepsilon_{M_c}]$ for some $\varepsilon_{M_c}>0$, hence in the interval $(d^*_+(c),d_{\max}-\varepsilon_{M_c}]$, for every $x\leq d_-^*(c)$. 
\end{itemize}
Hence we have shown that, given $x\in\Oc$, there exists a unique solution $y^*(x;c)$ to $L_1(x,\cdot;c)=0$ if and only if $c<\bar{c}_{+,g}$,
and it belongs to  the interval $(d^*_+(c),d_{\max}-\varepsilon_{M_c}]$.
Morever, Implicit Function Theorem ensures that $y^*(\cdot;c)$ is continuously differentiable and 
\beq\label{y**}
\frac{\ud}{\ud x} y^*(x;c)\ \ =\ \ - \ \frac{\frac{\partial L_1}{\partial x}(x,y^*(x;c))}{\frac{\partial L_1}{\partial y}(x,y^*(x;c))}\ \ =\ \ \frac{\psi(x)m'(x)(g_c(c,x)-\rho q_0^-)}{\psi(y^*(x))m'(y^*(x))(g_c(c,y^*(x))+\rho q_0^+)}.
\enq

Now consider the equation $L_2(x,y^*(x;c);c)=0$. We are going to show  existence and uniqueness of solutions to such equation in $\Oc$. This will complete the proof of existence and uniqueness of solutions for  \eqref{optcond0}, as, from what we have said  before, $x^*(c)$ solves the latter equation if and only if 
 $(x^*(c),y^*(x^*(c);c))$ solves \eqref{optcond0}.
We observe that:
\begin{itemize}
\item[-] If $c\leq \underline{c}_{-,g}$, then $g_c(c,\cdot)-\rho q_0^-\leq 0$ in $\Oc$; so, since $\varphi'(\cdot)/S'(\cdot)<0$ we have $L_2(\cdot,y^*(\cdot;c))<0$ in $\Oc$ and the solution does not exist.
\item[-] If $c>\underline{c}_{-,g}$, then we have the following facts:
\begin{itemize}
\item[1.] $L_2(\cdot;y^*(\cdot;c))<0$ in $(d^*_-(c),d_{\max})$, as $g_c(c,\cdot)-\rho q_0^-\leq 0$ therein and $\varphi'(\cdot)/S'(\cdot)<0$.
\item[2.] Using \eqref{y**} we compute
\beqs
\frac{\ud}{\ud x} L_2(x,y^*(x;c))\ \ =\ \ \frac{\psi(x)\varphi(y^*(x;c))-\psi(y^*(x;c))\varphi(x)}{\psi(y^*(x;c)}m'(x)(g_c(c,x)-\rho q_0^-).
\enqs
So taking into account that $y^*(x;c)>x$, the strict (opposite) monotonicity of $\varphi,\psi$, and that $g(c,\cdot)-\rho q_0^->0$ in $(d_{\min}, d^*_-(c))$, we see that $\frac{\ud}{\ud x} L_2(x,y^*(x;c))<0$ for $x\in(d_{\min}, d^*_-(c))$.
  \item[3.] Arguing as in proving \eqref{rrt}, we can prove that  there exists $\hat{x}\in (d_{\min},d^*_-(c))$ such that $L_2(\hat{x},y^*(\hat{x};c))<0$ and 
  \beqs
  L_2(x,y^*(x;c)) & \geq & \int_{\hat{x}}^{y^*(x;c)}   \varphi(\xi)(g_c(c,\xi)-\rho q_0^-)m'(\xi)\ud \xi\\
  &&-\frac{g_c(c,\hat{x})-\rho q_0^-}{\rho} \left(\frac{\varphi'(x)}{S'(x)}-\frac{\varphi'(\hat{x})}{S'(\hat{x})}\right).
  \enqs
  Since $y^*(x;c)\in (d^*_+(c), d_{\max}-\varepsilon_{M_c}]$  for every $x\in(d_{\min},d^*_-(c)]$, setting 
  $$K_0\ \ := \ \ \int_{\hat{x}}^{d_{\max}- \varepsilon_{M_c}}   \varphi(\xi)(g_c(c,\xi)-\rho q_0^-)m'(\xi)\ud \xi,$$ the latter inequality yields 
  \beqs
  L_2(x,y^*(x;c)) & \geq &K_0-\frac{g_c(c,\hat{x})-\rho q_0^-}{\rho} \left(\frac{\varphi'(x)}{S'(x)}-\frac{\varphi'(\hat{x})}{S'(\hat{x})}\right).
  \enqs
Now, since $\frac{\varphi'(x)}{S'(x)}\rightarrow -\infty$ as $x\rightarrow d_{\min}$ due to \eqref{psiphi'}, and since $g_c(c,\hat{x})-\rho q_0^->0$, we see that $L_2(x,y^*(x;c))\rightarrow \infty$ as $x\rightarrow d_{\min}$.
\end{itemize}
Combining these three fact we deduce that there exists a unique solution to the equation $L_2(\cdot;y^*(\cdot;c))=0$ and that it belongs to the interval $(d_{\min},d^*_-(c))$.\ep
\end{itemize}

Let us show now the last part of the claim. Consider $c$ as a variable in $L_1,L_2$ and consider the matrix
$$
\mathcal{M}(x^*(c),y^*(c);c) \ \ =\ \ \left(\begin{array}{cc} \frac{\partial L_1}{\partial x}(x^*(c),y^*(c);c) & \frac{\partial L_1}{\partial y}(x^*(c),y^*(c);c)\\
\frac{\partial L_2}{\partial x}(x^*(c),y^*(c);c)& \frac{\partial L_2}{\partial y}(x^*(c),y^*(c);c)\end{array}\right).
$$
Taking into account that $x^*(c)<d^*_-(c)$, $y^*(c)>d^*_-(c)$, and that $\psi,\varphi$ are respectively strictly increasing and strictly decreasing,  we see that the $\mathcal{M}^*(x^*(c),y^*(c);c)$  is actually non singular. More precisely $M:= \mbox{det}(\mathcal{M}(x^*(c),y^*(c);c))<0$ and
\beq\label{pppsd}
\mathcal{M}(x^*(c),y^*(c);c)^{-1} &=& \frac{1}{M}  \left(\begin{array}{cc} \frac{\partial L_2}{\partial y}(x^*(c),y^*(c);c) & -\frac{\partial L_1}{\partial y}(x^*(c),y^*(c);c)\\
-\frac{\partial L_2}{\partial x}(x^*(c),y^*(c);c)& \frac{\partial L_1}{\partial x}(x^*(c),y^*(c);c)\end{array}\right).
\enq
So, 
since $L_1(x^*(c),y^*(c);c)=0$, $L_1(x^*(c),y^*(c);c)=0$, we can  apply  Implicit Function Theorem which yields
\beq\label{x*y*}
\frac{\ud}{\ud c} \left(\!\!\!
\begin{array}{c}
x^*(c)\\ y^*(c)
\end{array}\!\!\!\right) \ \ =\ \ - \,\mathcal{M}(x^*(c),y^*(c);c)^{-1} \left(\begin{array}{c}
 \frac{\partial L_1}{\partial c}(x^*(c),y^*(c);c) \\ \frac{\partial L_2}{\partial c}(x^*(c),y^*(c);c)
\end{array}\right).
\enq
Since $g_{cc}=1$, we have
$$
\frac{\partial L_1}{\partial c}(x^*(c),y^*(c);c)\ \ =\ \ \int_{x^*(c)}^{y^*(c)} \psi(\xi)m'(\xi)\ud \xi, \ \ \ \ \frac{\partial L_2}{\partial c}(x^*(c),y^*(c);c)\ \ =\ \ \int_{x^*(c)}^{y^*(c)} \varphi(\xi)m'(\xi)\ud \xi. 
$$
So, from \eqref{x*y*}-\eqref{pppsd} we get
$$
\frac{\ud}{\ud c} x^*(c)\ \ =\ \ -\frac{1}{M}(g_c(c,y^*(c))+\rho q_0^+)m'(y^*(c))\int_{x^*(c)}^{y^*(c)}(\varphi(y^*(c))\psi(\xi)-\psi(y^*(c))\varphi(\xi))m'(\xi)\ud\xi
$$
$$
\frac{\ud}{\ud c} \,y^*(c)\ \ =\ \ -\frac{1}{M}(g_c(c,x^*(c))-\rho q_0^-)m'(x^*(c))\int_{x^*(c)}^{y^*(c)}(\varphi(x^*(c))\psi(\xi)-\psi(x^*(c))\varphi(\xi))m'(\xi)\ud\xi.
$$
Now, notice that 
$$M\  < \ 0, \ \ \ \ \ g_c(c,y^*(c))+\rho q_0^+\ < \ 0,  \ \ \ \  g_c(c,x^*(c))-\rho q_0^-\ >\ 0,$$
and that the functions 
$$q(\xi) \ : = \ \varphi(y^*(c))\psi(\xi)-\psi(y^*(c))\varphi(\xi), \ \ \ \ \ p(\xi) \ :=\ \varphi(x^*(c))\psi(\xi)-\psi(x^*(c))\varphi(\xi),$$
are both strictly increasing and verify, respectively $q(y^*(c))=0$ and $p(x^*(c))=0$. So we conclude from \eqref{x*y*}.
\ep\\\\
We are now ready to  characterize the optimal boundaries. 

\begin{Theorem}\label{teo:car}
Let $-\beta_0$ be strictly decreasing. We have $\underline{c}_-=\underline{c}_{-,g}$, $\bar{c}_+=\bar{c}_{+,g}$ and the optimal boundaries $\partial^\pm \Cc$ are characterized piecewise as follows. (Note that some of the three regions below where we split the characterization  may be empty.)
\begin{enumerate}
\item In the region 
$(\underline{c}_{-,g}, \bar{c}_{+,g})\times \Oc$, the optimal boundaries $\partial^\pm\Cc$ are identified by the functions $\hat{d}_\pm$  which are characterized as follows: given $c\in (\underline{c}_{-,g}, \bar{c}_{+,g})$ the couple $(\hat{d}_-(c),\hat{d}_+(c))\in \Oc\times\Oc$ is the unique solution  of the system of equations \eqref{optcond0} provided by Proposition \ref{prop:exun}.
\item In the region $(-\infty, \underline{c}_{-,g}]\times \Oc$ only $\partial^+\Cc$ (at most) exists and is identified in terms of the function $\hat{c}_+$ (note that  Corollary \ref{cor:fin} ensures $\hat{c}_+>-\infty)$), which is explicitly given by
$${\hat{c}_+}(d) \ \ = \ \ \rho \left[\beta(d)-\frac{\psi(d)}{\psi'(d)}\beta'(d)-{q_0^+}\right], \ \ \  d\ \leq \ \lim_{c\downarrow \underline{c}_{-,g}}\hat{d}_+(c), \ \ d\in\Oc.
$$
(For the definition of $ \lim_{c\downarrow \bar{c}_{-,g}}\hat{d}_+(c)$ when $(\underline{c}_{-,g}, \bar{c}_{+,g})$ is empty, recall that $\hat{d}_+(c)\equiv d_{\max}$ for $c\geq\bar{c}_+$.) 

\item In the region $[\bar{c}_{+,g},\infty)\times \Oc$ only $\partial^-\Cc$ (at most) exists and is identified in terms of the function $\hat{c}_-$ (note that  Corollary \ref{cor:fin} ensures $\hat{c}_-<\infty)$), which is explicitly given by
$${\hat{c}_-}(d) \ \ = \ \ \rho \left[\beta(d)-\frac{\varphi(d)}{\varphi'(d)}\beta'(d)+{q_0^-}\right], \ \ \ d\ \geq \ \lim_{c\uparrow \bar{c}_{+,g}}\hat{d}_-(c), \ \ d\in\Oc.
$$ 
(For the definition of $ \lim_{c\uparrow \bar{c}_{+,g}}\hat{d}_-(c)$ when $(\underline{c}_{-,g}, \bar{c}_{+,g})$ is empty, recall  that  $\hat{d}_-(c)\equiv d_{\min}$ for $c\leq\underline{c}_-$.) 
\end{enumerate}
Moreover:
\begin{itemize}
\item[(i)] The functions $\hat{c}_{\pm}:\mathcal{O}\rightarrow \R$ are   continuous and strictly increasing. 
\item[(ii)] $\hat{c}_+$ and $\hat{c}_-$ are of class $C^1$ except, at most, at the points $\lim_{c\downarrow \underline{c}_{-,g}}\hat{d}_+(c)$ and $\lim_{c\uparrow \bar{c}_{+,g}}\hat{d}_-(c)$, respectively (if they belong to $\Oc$).
\end{itemize}
\end{Theorem}
\noindent \textbf{Proof.}
1. First of all we notice that, by Proposition \ref{prop:stime}, we have $\underline{c}_{-,g}\leq \underline{c}_-$ and $\bar{c}_{+,g}\geq \bar{c}_+$. In the interval $(\underline{c}_-,\bar{c}_+)$, we have that the couple $(\hat{d}_-(c),\hat{d}_+(c))$ belongs to $\Oc\times\Oc$, and,  by Propositions \ref{prop:ex} and \ref{prop:exun}, it can  be identified as the unique solution  of the system of equations \eqref{optcond0}.  This shows claim 1, once we prove the claim $\underline{c}_{-,g}= \underline{c}_-$ and $\bar{c}_{+,g}= \bar{c}_+$, which is what we are going to prove now.

Assume by contradiction that $(\underline{c}_{-,g},\underline{c}_-]$ is nonempty. Then, for all $c\in (\underline{c}_{-,g},\underline{c}_-]$  we should have a unique solution $(d_-(c),d_+(c))\in\Oc\times\Oc$  to \eqref{optcond0} as provided by Proposition \ref{prop:exun}. By the monotonicity claim of Proposition \ref{prop:exun}, such a  solution should be such that $d_{\min}<d_-(c)<\lim_{\zeta\downarrow\underline{c}_-} \hat{d}_-(\zeta)=:d_0$. Now if $d_0>d_{\min}$, then, by definition of $\underline{c}_-$ we would have $\hat{c}_-\equiv \underline{c}_-$ in $(d_{\min},d_0)$ and we would have, by Proposition \ref{prop:ex}, more than one solution to \eqref{optcond0} at the level $\underline{c}_{-}$. But this contradicts Proposition \ref{prop:exun}. Therefore it should be $d_0=d_{\min}$, but  this would be  a contradiction to $d_{\min}<d_-(c)<d_0$. Hence,  it remains proved that $\underline{c}_{-}=\underline{c}_{-,g}$. The same argument applies to $\bar{c}_+$ and so the claim is proved.
\smallskip

2. The fact that only $\partial^+\Cc$ (at most) exists in the region $(-\infty,\underline{c}_{-,g}]$ is due to the definition of $\underline{c}_-$, to the equality $\underline{c}_{-,g}=\underline{c}_-$  and to the fact that, as shown in item 1,  $\lim_{\zeta\downarrow\underline{c}_-} \hat{d}_-(\zeta)=d_{\min}$. Then, due to Theorem \ref{prop:structure}, we have $B(c)=0$ for all $c\leq \underline{c}_{-,g}$. Hence, the optimality conditions  \eqref{smoothbase} and \eqref{opt2}  written at the points $(\hat{c}_+(d),d)\in\partial^+ \mathcal{C}$ with $d\in(d_{\min} , \lim_{c\downarrow \underline{c}_{-,g}}\hat{d}_+(c)]$ (notice that, due to Corollary \ref{cor:fin}, we actually have $\hat{c}_+:\Oc\rightarrow \R$) yield
  \begin{equation}\label{smoothbase3}
 \begin{cases}
 A'({\hat{c}_+}(d))\psi(d)+\frac{1}{\rho}{\hat{c}_+}(d)-\beta(d) \ \ = \ \ -{q_0^+},\\\\
A'({\hat{c}_+}(d))\psi'(d)-\beta'(d) \ \ =\ \  0.
 \end{cases}
 \end{equation}
 Multiplying the second equation in \eqref{smoothbase3} by $\psi/\psi'$ and subtracting it   to the first one, we get \eqref{hatcmr}. 
\smallskip

3. The same argument of item 2 applies symmetrically.

\medskip

Let us now show items (i) and (ii).

(i) We show the claim for $\hat{c}_+$, the proof of the claim regarding $\hat{c}_-$ is analogous. 

Since $\hat{d}_+$ is strictly increasing and continuous in the interval $(\underline{c}_{-,g},\bar{c}_{+,g})$ (when this is not empty), we see that $\hat{c}_+$ is the inverse of $\hat{d}_+$ in the interval $(\lim_{c\downarrow \underline{c}}\hat{d}_+(c),d_{\max})$ (when this is, correspondingly, nonempty) and is strictly increasing and continuous therein. So we must now show that $\hat{c}_+$ is strictly increasing and continuous in the interval $(d_{\min}, \lim_{c\downarrow \underline{c}}\hat{d}_+(c)]$ (when this is nonempty). Assume by contradiction that there exists a nonempty interval $(a,b)\subset (d_{\min}, \lim_{c\downarrow \underline{c}}\hat{d}_+(c)]$ where $\hat{c}_+\equiv c_0$. Then from the first equality in \eqref{smoothbase3} we should have
$$\beta(d)\ \ = \ \ A'(c_0)\psi(d)+\frac{1}{\rho}c_0 +q_0^+, \  \ \ \ d\in(a,b).$$
Since $\psi$ solves $\mathcal{L}\psi=0$, we then have that $\mathcal{L}\beta\equiv c_0+\rho q_0^+$ in $(a,b)$. On the other hand, from \eqref{alphabeta}, we see that it must be also $\mathcal{L}\beta=\beta_0$, so we should conclude that $\beta_0$ is constant in $(a,b)$, contradictiong the hypothesis. So, it has been proved that $\hat{c}_+$ is strictly increasing. 

Now we show that $\hat{c}_+$ is continuous. Indeed it is continuous in  the interval $(d_{\min}, \lim_{c\downarrow \underline{c}}\hat{d}_+(c)]$ due to item 2, and in the interval $(\lim_{c\downarrow \underline{c}}\hat{d}_+(c),d_{\max})$, due to item 1. It remains to prove that $\hat{c}_+$ is continuous at $\lim_{c\downarrow \underline{c}}\hat{d}_+(c)$ (when it belongs to $\Oc$). This comes just from the fact that $\hat{c}_+$ is right-continuous in general and, as we have seen just now, it is left-continuous at $\lim_{c\downarrow \underline{c}}\hat{d}_+(c)$.

(ii) It follows from the previuos claims and from Proposition \ref{prop:exun}.
\ep\\\\
We notice that $\underline{c}_{-,g}$, $\bar{c}_{+,g}$ are explicit.
So Theorem \ref{teo:car} actually provides a way to find, up to the (possibly numerical) solution of the system of equations \eqref{optcond0} for every $c\in (\underline{c}_{-,g},\bar{c}_{+,g})$, when this interval is not empty, the optimal boundaries $\partial^\pm\Cc$. Then the functions $A,B$ individuating the value function in the continuation region can be retrieved by Theorem \ref{prop:structure}: 
\begin{itemize} 
\item[-] If  $(\underline{c}_{-,g},\bar{c}_{+,g})\neq \emptyset$,  then
$A,B$ can be computed in the interval $(\underline{c}_{-,g},\bar{c}_{+,g})$ by integrating \eqref{eds} with boundary conditions $A(\bar{c}_{+,g})=0$ and $B(\underline{c}_{-,g})=0$, and, respectively in the intervals $(\underline{c}_+,\underline{c}_{-,g}]$ and $[\bar{c}_+,\bar{c}_{+,g})$ (when they are nonempty), by the equalities
$$
A(c) \ \ =\ \ [\lim_{c\downarrow \underline{c}_{-,g}}A(c)] -\int_c^{\underline{c}_{-,g}} \frac{\beta'(\hat{d}_+(\xi))}{\psi'(\hat{d}_+(\xi))}\ud\xi, \ \ \ \ c\in (\underline{c}_+,\underline{c}_{-,g}]; $$
$$B(c) \ \ =\ \  [\lim_{c\uparrow \bar{c}_{+,g}}B(c)] +\int_{\bar{c}_{+,g}}^c \frac{\beta'(\hat{d}_+(\xi))}{\varphi'(\hat{d}_+(\xi))}\ud\xi, \ \ \ \ c\in [\bar{c}_{+,g}, \bar{c}_-).
$$ 
\item[-] If $(\underline{c}_{-,g},\bar{c}_{+,g})= \emptyset$, then 
\begin{equation}\label{AA}
A(c) \ \ = \ \ -\int_c^{\bar{c}_{+,g}} \frac{\beta'(\hat{d}_+(\xi))}{\psi'(\hat{d}_+(\xi))}\ud\xi, \ \ \ \ c\in(\underline{c}_+,\bar{c}_{+,g}),
\end{equation}
\begin{equation}\label{BB}
B(c) \ \ = \ \ \int_{\underline{c}_{-,g}}^c \frac{\beta'(\hat{d}_-(\xi))}{\varphi'(\hat{d}_+(\xi))}\ud\xi, \ \ \ \ c\in(\underline{c}_{-,g},\bar{c}_-).
\end{equation}
\end{itemize}

Then $z_\pm$ can be obtained by \eqref{expz}.

\subsection{Quadratic cost and irreversibility}\label{sec:irr}
 In this subsection we consider we further particularize to the irreversible investment  case. Even if it is, rigorously speaking, out of our setting, nonetheless it can be formally seen as corresponding to take  $q_0^-=\infty$. The upper boundary in this case is clearly $\hat{c}_-\equiv\infty$, or, in other terms, it disappears. Hence, from Theorem \ref{teo:car}, we immediately get the following. 
 \begin{Corollary}\label{sol}
Let $q_0^-=\infty$,  and let the assumptions of Theorem \ref{teo:car} hold true. Then the functions $\hat{c}_{\pm},A,B, z_{{\pm}}$ of Theorem \ref{prop:structure} are  determined as follows: 
\begin{itemize}
\item[(a)] 
The upper optimal boundary is $\hat{c}_-\equiv\infty$, and 
lower boundary function ${\hat{c}_+}$ is  explicitly given by
\begin{equation}\label{hatcmr}
{\hat{c}_+}(d) \; = \; \rho \left[\beta(d)-\frac{\psi(d)}{\psi'(d)}\beta'(d)-{q_0^+}\right], \ \ \  d\in\mathcal{O}.
\end{equation}
{In particular ${\hat{c}_+}\in C^1(\mathcal{O};\mathbb{R})$.}
\item[(b)] $B\equiv 0$, and the function $A$ is given by
\begin{equation*}
A(c) \; = \; -\int_c^{\bar{c}_{+,g}} \frac{\beta'(\hat{d}_+(\xi))}{\psi'(\hat{d}_+(\xi))}\ud\xi, \ \ \ \ c\in(\underline{c}_+,\bar{c}_{+,g}),
\end{equation*}
\item[(c)]  The function $z_-$ is whatever function (it does not play a role, as $\hat{c}_- \equiv \infty$ implies $\mathcal{A}^-=\emptyset$), while the function $z_{+}$ is 
\beq\label{zz}
z_+(d) & = & A({\hat{c}_+}(d))\psi(d)+{{\hat{V}}}({\hat{c}_+}(d),d)+{q_0^+}{\hat{c}_+}(d), \ \  \ \ d\in\Oc,
\enq
with $\hat{V}$ given in \eqref{hatVqua}.
\end{itemize}
\end{Corollary}

\medskip
We end this paper by a simple and explicit illustration of our Corollary \ref{sol} to the case when the demand is modeled as a geometric Brownian motion:
\beqs 
 \ud D_t &=& \mu D_t{{\ud}t}+\sigma D_t\ud W_t, \ \ \ \mu\in\mathbb{R},\ \sigma>0,
 \enqs
 with initial datum $d>0$. In this case $\mathcal{O}=(0,\infty)$.
 Moreover, assume that
 \beqs
 g(c,d) &=& \frac{1}{2}(c-d)^2,
 \enqs
 and, according to  \eqref{rhoK1}, assume that
\begin{equation}\label{ass:rho2}
 \rho\ \ >\ \ [2\mu+\sigma^2]^+.
\end{equation}
 Then ${\hat{V}}$ is the quadratic function equal to
\beqs
{\hat{V}}(c,d) & = & \frac{1}{2}\Big(\frac{1}{\rho-2\mu-\sigma^2}d^2-\frac{2}{\rho-\mu}dc+\frac{1}{\rho}c^2\Big).
\enqs
The increasing fundamental solution to
 \beqs \label{A}
 [\mathcal{L}\phi](d) \ \  := \ \  \rho \phi-\mu d\phi'-\frac{1}{2}\sigma^2 d^2\phi'' &=& 0,
 \enqs
is given by
 \beqs
 \psi(d) &=& d^{m},
 \enqs
 where $m$ is the positive  root of the equation $\rho-\mu m-\frac{1}{2}\sigma^2 m(m-1)$ $=$ $0$,  and explicitly given by
 \beqs
 m &=& - \frac{\mu}{\sigma^2} + \frac{1}{2} + \sqrt{ \Big( - \frac{\mu}{\sigma^2} + \frac{1}{2}  \Big)^2 + \frac{ 2\rho}{\sigma^2} } 
 \enqs
 Notice that  $m$ $>$ $2$ by \eqref{ass:rho2}.  
 From Corollary \ref{sol},  the value function $v$ has the explicit form 
 $$v(c,d) \ \ = \ \ \begin{cases}A(c)d^m+{{\hat{V}}}(c,d), \ \ \ \ \ \ \ \ \mbox{if} \ c>{\hat{c}_+}(d),\\
 -{q_0^+}c+z(d), \ \ \ \ \ \ \ \ \ \ \ \ \, \ \ \mbox{if} \ c\leq {\hat{c}_+}(d).
 \end{cases}
 $$
 where the functions $A,{\hat{c}_+},z$ are 
 \beqs
{\hat{c}_+}(d) &= & a d-b, \ \ \ \ d>0, \\
A(c) &=&  - \frac{a^{m-1}}{m(m-2)}\frac{1}{\rho-\mu}(c+b)^{2-m}, \ \ \ \ \ c>-b, \\
z_+(d) &=& A(a d-b)d^m+{{\hat{V}}}(a d-b,d)+{q_0^+}(\alpha d-b), \ \ \ \ d>0,
\enqs
with
\beqs
a \; = \; \frac{m-1}{m}\frac{\rho}{\rho-\mu}, & &  b=\rho {q_0^+}.
\enqs 
 
 \appendix
 
\section{Appendix}

\noindent \textbf{Proof of Proposition \ref{exist}.}
\emph{Existence.} Let $(c,d)\in\mathcal{S}$ and take a sequence $(I^n)_{n\in\N}\subset\mathcal{I}$ s.t. $G(c,d;I^n)\rightarrow v(c,d)$. Assume, without loss of generality, that $G(c,d;I^n)\leq v(c,d)+1$ for all $n\geq 0$ and set   $\kappa:=\mbox{min}\{q_0^+,q_0^-\}>0$. Then, taking into account that $g\geq 0$, that $I_{0_-}^{n,+}=I_{0^-}^{n,-}=0$ for all $n\geq 0$, and integrating by parts, we get
\beqs
v(c,d)+1&\geq &\kappa \ \E \int_0^\infty e^{-\rho t}\left(\ud I_t^{n,+}+\ud I_t^{n,-}\right) \\
&= & \kappa\ \E \left[\int_0^\infty e^{-\rho t}\left(I_t^{n,+}+I_t^{n,-}\right)\ud t +[e^{-\rho t} (I_t^{n,+}+I_t^{n,-})]^\infty_{0^-}\right]\\
&\geq & \kappa\ \E \left[\int_0^\infty e^{-\rho t}\left(I_t^{n,+}+I_t^{n,-}\right)\ud t\right].
\enqs
So, the sequence $(I^n)_{n\in\N}$ is bounded in the space $L^1(\Omega\times \R;\P\times e^{-\rho t}\ud t)$. Thus, by a theorem of Koml\'os, there exists a subsequence (relabeled and still denoted by $(I^n)_{n\in\mathbb{N}}$) and a pair of measurable processes $\tilde{I}^+,\tilde{I}^-$ such that the Ces\`aro sequences of processes
\beq\label{Ipm}
\left(\tilde{I}^{n,\pm} \ :=\ \frac{1}{n}\sum_{j=1}^n I^{n,\pm}\right)\ \subset \ \mathcal{I} \ \ \ \mbox{converge} \ \ \ (\P\times e^{-\rho t}\ud t)-\mbox{a.e. to} \ \tilde{I}^{\pm}.
\enq
Define $\tilde{I}^n:=\tilde{I}^{n,+}-\tilde{I}^{n,-}$. Then, from \eqref{Ipm}, we have the convergence
\beq
\tilde{I}^{n} \ \ \longrightarrow \ \ \tilde{I} \ \ \ \ \ (\P\times e^{-\rho t}\ud t)-\mbox{a.e.}.
\enq
By convexity of $G$ w.r.t. the control argument $I$, we have that also $(\tilde{I}^n)_{n\in\N}$ is a minimizing sequence, i.e. $G(c,d;\tilde{I}^n)\rightarrow v(c,d)$.
On the other hand, arguing as in Lemmata 4.5--4.7 of \cite{KS1}, we can see that 
$\tilde{I}^+$ and $\tilde{I^-}$ admit modifications - which we still denote by   $\tilde{I}^+$ and $\tilde{I}^-$ - right-continuous, nondecreasing, and $\F$-adapted.
Hence, there is also a modification of $\tilde{I}$ - which we still denote by $\tilde{I}$ - belonging to $\mathcal{I}$. Now
Fatou's Lemma yields
\beq
G(c,d;\tilde{I})&\leq&\liminf_{n\rightarrow\infty}G(c,d;\tilde{I}^n)\ \ =\ \ v(c,d),
\enq
so $\tilde{I}$ is an optimal control starting from $(c,d)$.
\smallskip

\emph{Uniqueness.} Let $(c,d)\in\mathcal{S}$, and let  $I^{1}\in\mathcal{I}$, $I^{2}\in\mathcal{I}$ be two optimal controls starting from $(c,d)$. Define $\bar{I}:=\frac{1}{2}I^1+\frac{1}{2}I^2$. By linearity of the state equation \eqref{GGG} we then have  $C^{c,\bar{I}}=\frac{1}{2}C^{c,I^1}+\frac{1}{2}C^{c,I^2}$. Thus,  since $g(\cdot,d)$ is convex,
\beqs
0&\leq &G(c,d;\bar{I})-v(c,d)\ \
=\ \ G(c,d;\bar{I})-\frac{1}{2}G(c,d;I^1)-\frac{1}{2}G(c,d;I^2)\\&=&\E
\left[ \int_0^\infty e^{-\rho t} 
\left(
g\Big(\frac{1}{2}C^{c,I^1}+\frac{1}{2}C^{c,I^2},D_t^d\Big)-\frac{1}{2}g(C_t^{c,I^1},D_t^d)-\frac{1}{2}g(C_t^{c,I_2},D_t^d)
\right)
\right]\ \ \leq \ \ 0 .
\enqs
So, the inequalities above are indeed equalities and, still due to convexity of $g(\cdot,d)$, we must have
\beqs
g(C^{c,\bar{I}},D_t^d)-\frac{1}{2}g(C_t^{c,I^1},D_t^d)-\frac{1}{2}g(C_t^{c,I_2},D_t^d) &=& 0, \ \ \ \P-\mbox{a.s.}, \ \mbox{for a.e.} \ t\in\mathbb{R}. 
\enqs
Now the assumption of strict convexity of $g(\cdot,d)$ implies
$C^{c,I^1}=C^{c,I^2}$, $\P-\mbox{a.s.}$, for a.e. $t\in\mathbb{R}$, from which we derive $I^1=I^2$,  $\P-\mbox{a.s.}$, for a.e. $t\in\mathbb{R}$. So, due to right-continuity, $I^1$ and $I^2$ are indistinguishable.
\medskip\ep

\begin{Lemma}
Let $(c,d)\in\mathcal{S}$ and denote by $v_c^+(c,d)$, $v_c^-(c,d)$, respectively, the right- and left-derivative of $v$ w.r.t. $c$ at $(c,d)$ (their existence being guaranteed by convexity of $v(\cdot,d))$.
Then
\beq\label{lrder}
v_c^+(c,d)\ \ \leq \ \ J(c,d;\sigma, \tau^*), \ \ \ \forall \sigma\in\mathcal{T}; \ \ \ \ \  v_c^-(c,d)\ \ \geq \ \ J(c,d;\sigma^*, \tau), \ \ \ \forall \tau\in\mathcal{T}. \
\enq 
\end{Lemma}
\noindent \textbf{Proof.} Let us show the first inequality.
Let $(c,d)\in\Sc$ and let $I^*=(I^{*,+},I^{*,-})\in\mathcal{I}$ be an optimal control for $(c,d)$. Let $\varepsilon>0$ and set
\beqs
\tau^*\ :=\ \inf\{t\geq 0 \ | \ I_t^{*,+}>0\},  \ \ \ \ \ \   \tau_\varepsilon\ :=\ \inf\{t\geq 0 \ | \ I_t^{*,+}\geq \varepsilon\}.
\enqs
Moreover, given $\sigma\in\Tc$, set 
\beqs
I^\varepsilon &:=&\begin{cases}
-I_t^{*,-}, \ \ \ \ \mbox{if} \ \ 0\leq t<\sigma\wedge \tau_\varepsilon,\\
I^{*}_t-\varepsilon, \ \ \ \ \mbox{if}\ \ t\geq \sigma\wedge\tau_\varepsilon.
\end{cases}
\enqs
We can write
\beqs
G(c+\varepsilon,d;I^\varepsilon) \ \ = \ \ \E\bigg[\int_0^{\sigma\wedge\tau^*}e^{-\rho t} g(c+\varepsilon-I_t^{*,-},D_t^d)\ud t \\+\int_{\sigma\wedge\tau^*}^{\sigma\wedge \tau_\varepsilon} e^{-\rho t} g(c+\varepsilon-I_t^{*,-},D_t^d)\ud t
+\int_{\sigma\wedge \tau_\varepsilon}^\infty e^{-\rho t} g(c-I_t^{*},D_t^d)\ud t \\+\mathbf{1}_{\{\tau_\varepsilon\leq \sigma\}}\big(e^{-\rho\tau_\varepsilon}q_0^+(I_{\tau_\varepsilon}^{*,+}-\varepsilon)+\int_{\tau_\varepsilon^+}^\infty e^{-\rho t} q_0^+\ud I_t^{*,+}+\int_0^\infty e^{-\rho t}q_0^-\ud I_t^{*,-}\big)\\
+\mathbf{1}_{\{\tau^*\leq\sigma<\tau_\varepsilon\}}\big(e^{-\rho \sigma} q_0^-(\varepsilon-I_\sigma^{*,+})+\int_{\sigma^+}^\infty e^{-\rho t} q_0^+\ud I_t^{*,+}+\int_0^\infty e^{-\rho t}q_0^-\ud I_t^{*,-}\big)\\
+\mathbf{1}_{\{\sigma<\tau^*\}}\big(e^{-\rho \sigma} q_0^-\varepsilon+\int_{\tau^*}^\infty e^{-\rho t} q_0^+\ud I_t^{*,+}+\int_0^\infty e^{-\rho t}q_0^-\ud I_t^{*,-}\big)\bigg],
\enqs
and
\beqs
G(c,d;I^*)  = \ \ \E\bigg[\int_0^{\sigma\wedge\tau^*}e^{-\rho t} g(c-I_t^{*,-},D_t^d)\ud t +\int_{\sigma\wedge\tau^*}^{\sigma\wedge \tau_\varepsilon} e^{-\rho t} g(c+I_t^{*},D_t^d)\ud t\\
+\int_{\sigma\wedge \tau_\varepsilon}^\infty e^{-\rho t} g(c+I_t^{*,-},D_t^d)\ud t \\
+\mathbf{1}_{\{\tau_\varepsilon\leq \sigma\}}\Big(\int_{\tau^*}^{\tau_\varepsilon^-}e^{-\rho t} q_0^+\ud I_t^{*,+}+e^{-\rho\tau_\varepsilon}q_0^+(I_{\tau_\varepsilon}^{*,+}-I_{\tau_\varepsilon^-}^{*,+})+\int_{\tau_\varepsilon^+}^\infty e^{-\rho t} q_0^+\ud I_t^{*,+}+\int_0^\infty e^{-\rho t}q_0^-\ud I_t^{*,-}\Big)\\
+\mathbf{1}_{\{\tau^*\leq\sigma<\tau_\varepsilon\}}\Big(\int_{\tau^*}^{\sigma^-}e^{-\rho t} q_0^+\ud I_t^{*,+}+e^{-\rho \sigma} q_0^-(I_\sigma^{*,+}-I_{\sigma^-}^{*,+})+\int_{\sigma^+}^\infty e^{-\rho t} q_0^+\ud I_t^{*,+}+\int_0^\infty e^{-\rho t}q_0^-\ud I_t^{*,-}\Big)\\
+\mathbf{1}_{\{\sigma<\tau^*\}}\Big(\int_{\tau^*}^\infty e^{-\rho t} q_0^+\ud I_t^{*,+}+\int_0^\infty e^{-\rho t}q_0^-\ud I_t^{*,-}\Big)\bigg].
\enqs
Subtracting we get
\beqs
v(c+\varepsilon,d) -v(c,d)\ \ \leq  \ \ \E\bigg[\int_0^{\sigma\wedge\tau^*}e^{-\rho t} \big(g(c+\varepsilon-I_t^{*,-},D_t^d)- g(c-I_t^{*,-},D_t^d)\big)\ud t \\
+\int_{\sigma\wedge\tau^*}^{\sigma\wedge \tau_\varepsilon} e^{-\rho t} \big(g(c+\varepsilon-I_t^{*,-},D_t^d)- g(c+I_t^{*,+}-I_t^*,D_t^d)\big)\ud t\\
+\mathbf{1}_{\{\tau_\varepsilon\leq \sigma\}}\Big(e^{-\rho\tau_\varepsilon}q_0^+(I_{\tau_\varepsilon^-}^{*,+}-\varepsilon)-\int_{\tau^*}^{\tau_\varepsilon^-}e^{-\rho t} q_0^+\ud I_t^{*,+}\Big)\\
+\mathbf{1}_{\{\tau^*\leq\sigma<\tau_\varepsilon\}}\Big(e^{-\rho \sigma} q_0^-(I_{\sigma^-}^{*,+}-\varepsilon)-\int_{\tau^*}^{\sigma^-} e^{-\rho t} q_0^+\ud I_t^{*,+}\Big)-\mathbf{1}_{\{\sigma<\tau^*\}}e^{-\rho \sigma} q_0^-\varepsilon\bigg].
\enqs
Using convexity of $g(\cdot,d)$ we can estimate from above the first two terms in the expectation  above respectively with
\beqs
 \varepsilon\ \int_0^{\sigma\wedge\tau^*}e^{-\rho t} g_c(c-I_t^{*,-},D_t^d)\ud t,\ \ \ \ \ \ L_1(\varepsilon)\ :=\ \int_{\sigma\wedge\tau^*}^{\sigma\wedge\tau_\varepsilon}e^{-\rho t} (\varepsilon-I_t^{*,+})g_c(c+\varepsilon,D_t^d)\ud t,
\enqs
while the third term can be rearranged as
\beqs
-\varepsilon q_0^+ e^{-\rho\tau^*}\mathbf{1}_{\{\tau^*<\sigma\}} + L_2(\varepsilon) +L_3(\varepsilon),
\enqs
where
\beqs
L_2(\varepsilon)\ :=\ \varepsilon q_0^+[e^{-\rho\tau^*}\mathbf{1}_{\{\tau^*<\sigma\}}-e^{-\rho \tau_\varepsilon}\mathbf{1}_{\{\tau_\varepsilon\leq \sigma\}}], \ \ \ \  L_3(\varepsilon)\ := \mathbf{1}_{\{\tau_\varepsilon\leq \sigma\}}\Big(e^{-\rho\tau_\varepsilon}I^{*,+}_{\tau_\varepsilon^-}-\int_{\tau^*}^{\tau_\varepsilon^-}e^{-\rho t}\ud I_t^{*,+}\Big).
\enqs
Setting also 
\beqs
L_4(\varepsilon)\ :=\ \mathbf{1}_{\{\tau^*\leq\sigma<\tau_\varepsilon\}}\Big(e^{-\rho \sigma} q_0^-(I_{\sigma^-}^{*,+}-\varepsilon)-\int_{\tau^*}^{\sigma^-} e^{-\rho t} q_0^+\ud I_t^{*,+}\Big)
\enqs
we can write
\beqs
\frac{v(c+\varepsilon,d) -v(c,d)}{\varepsilon} &\leq & J(c,d;\sigma,\tau^*)+\frac{1}{\varepsilon}\sum_{j=1}^{4} L_j(\varepsilon).
\enqs
Using estimates like the ones used in \cite[Lemma 4.3]{KW}, one can see that, for each $j=1,...,4$,
$\frac{L_j(\varepsilon)}{\varepsilon}\rightarrow 0$ when $\varepsilon\to 0$, which gives the first inequality in \eqref{lrder}. The second inequality can be obtained in a similar way.
\ep\medskip\\
\noindent \textbf{Proof of Proposition \ref{Prop:dynkin}.} Let $v^+_c(c,d)$ and $v_c^-(c,d)$ be, respectively, the left and the right derivative of $v$ w.r.t. $c$ at $(c,d)$, which exist due to convexity of $v(\cdot,d)$ and verify $v_c^+(c,d)\geq v_c^-(c,d)$.
Then, considering \eqref{lrder}, we get 
\beqs
v_c^-(c,d)\ \ \leq \ \ v_c^+(c,d)\ \ \leq\ \  J(c,d;\sigma^*,\tau^*)\ \ \leq \ \ v_c^-(c,d).
\enqs
So the inequalities above are indeed equalities and hence it follows that $v_c(c,d)$ exists and is equal to $J(c,d;\sigma^*,\tau^*)$. Then, still using \eqref{lrder}, we get
\beqs
J(c,d;\sigma^*,\tau)\ \leq \ v^-_c(c,d)\ = \ J(c,d;\sigma^*,\tau^*) \ =\ v_c^+(c,d) \ \leq \ J(c,d;\sigma,\tau^*),\ \ \ \ \forall \sigma\in\mathcal{T}, \ \forall \tau\in\mathcal{T}.
\enqs
This shows both the claims.
\ep

 \begin{Proposition}[It\^o's Formula]\label{prop:ITO}
 Let $\varphi\in C^{1,2}(\mathcal{S};\R)$, $(c,d)\in\mathcal{S}$, $I\in\mathcal{I}$, and let $\tau$ be a bounded stopping time such that $(C_t^{c,I}, D_t^d)_{t\in[0,\tau]}$ is contained in a compact subset of $\mathcal{S}$. Then the following change of variable's formula holds:
 \beq
\varphi(c,d) &=& \mathbb{E}\Big[e^{-\rho\tau}\varphi(C^{c,I}_{\tau},D^d_{\tau})\Big] 
\;   +  \;  \mathbb{E}\Big[\int_0^{\tau}e^{-\rho t}[\mathcal{L}\varphi(C_t^{c.I},\cdot)](D_t^d){{\ud}t}\Big]\nonumber\\
& &  -   \; \mathbb{E}\Big[\int_0^{\tau}e^{-\rho t} \varphi_c(C_t^{c,I},D_t^d)\ud I_t\Big] \;  \nonumber \\
& & - \; \mathbb{E}\Big[\sum_{0\leq t\leq \tau}e^{-\rho t}(\varphi(C_t^{c,I},D_t^d)-\varphi(C_{t^-}^{c,I},D_t^d)-\varphi_c(C_t^*,D_t^d) \Delta C_t^{c,I})\Big], \nonumber 
\enq
 \end{Proposition}
\noindent  \textbf{Proof.} Theorem 33 (p.\,81) in \cite{Pr} provides the desired formula for functions which are continuously twice differentiable when $\tau$ is constant. The extension to the case of $\tau$ stopping time for the latter class of functions is standard. To get the formula for functions belonging to  $C^{1,2}(\mathcal{S};\R)$, one can argue using mollifiers as follows.
Take a sequence of mollifiers $(\xi_n)_{n\in\mathbb{N}}$ and consider the convolution $\varphi_n:=\xi_n\ast \varphi$.  Then $\varphi_n$ is continuously twice differentiable for each $n$, so the formula applies to the sequence $(\varphi_n)_{n\in\mathbb{N}}$. Moreover all the derivatives of $\varphi_n$ involved in the formula converge locally uniformly  to the corresponding derivatives of $v$ (which exist, as the formula involves only derivatives which are defined in the class   
$C^{1,2}(\mathcal{S};\R)$. Hence, the claim follows by uniform convergence since $(C_t^{c,I}, D_t^d)_{t\in[0,\tau]}$ is contained in a compact subset of $\mathcal{S}$. \ep
 
%
%
 \vspace{3mm}
 
\noindent {\small \textbf{Acknowledgement.} The authors would like to thank Ren\'e A\"id and Bertrand Villeneuve for useful discussions. 
Salvatore Federico also thanks  Giorgio Ferrari for comments  on the literature and for suggesting  the possible extension of the model described in Remark \ref{giorgio}.}

\end{document}